\documentclass[reqno]{amsart}

\usepackage{geometry}
\usepackage[english]{babel}
\usepackage{latexsym, mathrsfs, verbatim}
\usepackage{amsfonts, amsthm, amsmath, amssymb, bm}
\usepackage{amsaddr}

\usepackage{verbatim}
\usepackage{mathrsfs,amsmath,enumerate}
\usepackage{pgf}
\usepackage{comment}

\newtheorem{teor}{Theorem}[section]

\newtheorem{prop}[teor]{Proposition}
\newtheorem{lemma}[teor]{Lemma}

\theoremstyle{definition}

\newcommand{\F}{\mathbb{F}}
\newcommand{\N}{\mathbb{N}}
\newcommand{\Z}{\mathbb{Z}}

\newcommand{\C}{\mathbb{C}} 
\newcommand{\ZZ}{\mathbf{Z}}

\newcommand{\Aut}{\mathrm{Aut}}

\newcommand{\PSL}{\mathrm{PSL}}
\newcommand{\SL}{\mathrm{SL}}
\newcommand{\GU}{\mathrm{GU}}
\newcommand{\SU}{\mathrm{SU}}
\newcommand{\PSU}{\mathrm{PSU}}
\newcommand{\PGU}{\mathrm{PGU}}

\newcommand{\PGL}{\mathrm{PGL}}
\newcommand{\GL}{\mathrm{GL}}

\newcommand{\Sz}{\mathrm{Sz}}

\newcommand{\soc}{\mathrm{soc}}

\newcommand{\Irr}{\mathrm{Irr}}

\def\equad{\quad \text{ and } \quad}

\numberwithin{equation}{section}
\numberwithin{table}{section}

\usepackage[table,dvipsnames]{xcolor}

\usepackage{hyperref}
\hypersetup{ 
 colorlinks =true,
 linkcolor=Sepia,
 citecolor=Plum,
 urlcolor=olive,
 }

\begin{document}

\title{Unisingular representations of rank 1 finite simple groups of Lie type}

\author[Marco Antonio Pellegrini]{Marco Antonio Pellegrini}
\address{Dipartimento di Matematica e Fisica, Universit\`a Cattolica del Sacro Cuore,\\
Via della Garzetta 48, 25133 Brescia, Italy}
\email{marcoantonio.pellegrini@unicatt.it}

\author[Lorenzo Schena]{Lorenzo Schena}
\address{Lehrstuhl für Algebra und Zahlentheorie,
RWTH Aachen University, \\
Pontdriesch 14-16, 52062 Aachen, Germany}
\email{lorenzo.schena@rwth-aachen.de}

\begin{abstract}
A representation $\Phi: G \to \GL_n(\F)$ of a finite group $G$ is called unisingular if the matrix $\Phi(g)$ admits $1$ as an eigenvalue for any $g\in G$.
In this paper, we determine all the complex irreducible unisingular representations of the finite simple groups of Lie type of rank $1$ and of the almost simple sporadic groups.
\end{abstract}

\maketitle

\section{Introduction}

Several questions on finite groups require the knowledge of the eigenvalues of the matrices that arise from their irreducible representations (for instance, see \cite{DPZ,DZ}). In particular, it can be of particular interest to know when,
given an irreducible representation $\Phi$  of a finite group $G$ and an element $g \in G$, the matrix $\Phi(g)$ admits the eigenvalue $1$.

For instance, when in 1935 Hans Zassenhaus \cite{Zas} classified the near-fields, he also gave a complete description of
those finite groups which have complex representations $\Phi$ such that $\Phi(g)$ does not have the eigenvalue $1$ for any nontrivial element $g \in G$.
As an opposite problem, one can study finite groups whose matrix representations always have the eigenvalue $1$.
An irreducible representation $\Phi$  of  $G$ (and the associated character) is called \emph{unisingular} if $\Phi(g)$ has the eigenvalue $1$ for every $g\in G$.

The term unisingular was introduced
by László Babai and Aner Shalev: in \cite{BS} a finite group is said to be \emph{unisingular in characteristic $p$} if every its element admits a non-zero invariant vector in every irreducible representation in characteristic $p$. 
Robert Guralnick and Pham Huu  Tiep classified in \cite{GT} the finite simple groups of Lie type of characteristic $p$ that are unisingular in characteristic $p$.

More in general, John Cullinan and Alexandre Zalesski  in \cite{CZ} used the term unisingular to denote any representation $\Phi: G \to \GL_n(\F)$ of a group $G$ such that for any $g\in G$ the matrix $\Phi(g)$ admits $1$ as an eigenvalue.
We shall use the term unisingular with this meaning, but focusing our attention to complex irreducible representations.
We shall also speak of unisingular character with the  meaning that the associated (complex) representations are unisingular.

Very few is known about complex unisingular representations of finite simple groups.
In addition to the complete analysis for the alternating 
(and symmetric) groups provided in \cite{Am1,Am}, only some finite simple groups of Lie type have been studied, mainly by Zalesski in 
\cite{Z88,Z90,Z09} (see also \cite{Z25} for a recent survey on this topic). Hence, determining all the irreducible unisingular representations of a finite group is still a challenging open problem.

In this paper, we determine the complex unisingular representations of the finite simple groups of Lie type of rank $1$. We recall that these groups are  the linear groups $\PSL_2(q)$, the unitary groups $\PSU_3(q)$, the Suzuki groups $\Sz(2^{2m+1})$, and the small Ree groups  ${}^2\mathrm{G}_2 (3^{2m+1})$. Our main result is the following.

\begin{teor}\label{thm:rk1}
Let $G$ be a rank $1$ finite simple group of Lie type.
An irreducible character $\chi$ of $G$ is unisingular if and only if none of the following cases occurs:
\begin{enumerate}[$(1)$]
\item $G=\PSL_2(q)$,  $\chi(1)=\frac{q-1}{2}$ and  $q\geq 7$ is an odd prime such that $q\equiv 3 \pmod 4$;

\item $G=\PSL_2(q)$,  $\chi(1)=q-1$ and $q\geq 5$ is an odd prime;

\item $G=\PSL_2(q)$, $\chi(1)=q$ and $q\geq 4$ is even;

\item $G=\PSU_3(q)$, $\chi(1)=q^2-q$
and $q\geq 3$ is an odd prime if $q\equiv 2 \pmod 3$;

\item $G=\PSU_3(q)$,  $\chi(1)=q^2-q+1$, $q\geq 3$ is an odd prime and $\chi$ is not rational-valued.

\end{enumerate}
\end{teor}

We also consider the groups $\PGL_2(q)$ and $\PGU_3(q)$
(see  Propositions \ref{pPGL2} and \ref{pPGU3}).
In the last section, we classify the complex unisingular representations of the almost simple sporadic groups (see Proposition \ref{spor}).

\section{Preliminary results}

Let $\Phi: G \to \GL_n(\C)$ be a representation of a finite group $G$.
If we know the power maps of any element $g\in G$, we can compute the algebraic multiplicity of any eigenvalue of the matrix $\Phi(g)$ from the character table of $G$ by applying the following.

\begin{prop}
Let $G$ be a group, $\Phi$ be a representation of $G$, $\chi$ be the character afforded by $\Phi$, $g\in G$, $m$ be the order of $g$ and $\zeta_m\in \C^*$ be a primitive $m$-root of $1$.
Then, for any $i\in \N$, the algebraic multiplicity of $\zeta_m^i$ as eigenvalue of $\Phi(g)$ is
$\dfrac{1}{m}\sum\limits_{j=0}^{m-1}\chi(g^j)\zeta_m^{-ij}$.
\end{prop}

In particular, if $\chi$ is the character of $\Phi$ and $g$ is an element of $G$ of order $|g|=m$, the multiplicity of the eigenvalue $1$ in $\Phi(g)$ is given by the formula
\begin{equation} \label{MolteplicityOf1}
		M_\chi(g)=\frac{1}{m}\sum_{j=0}^{m-1}\chi(g^j).
\end{equation}
By Frobenius' reciprocity, we have $M_\chi(g)= \left[\chi_{|\langle g \rangle},1_{\langle g \rangle}\right]_{\langle g \rangle} =
\left[\chi,(1_{\langle g \rangle})^G\right]_G$.

It is immediate to see that the principal character is always unisingular and it is the only unisingular linear character.
Now, we discuss the connection between unisingular characters and the center of the group.
Recall that the \emph{quasikernel} of a character $\chi$  of a group $G$ is the set
$\ZZ(\chi) = \{g\in G : |\chi(g)|= \chi(1) \}$.

\begin{lemma}\label{ZK}
If an irreducible character $\chi$ of a finite group $G$ is unisingular, then $\ZZ(\chi)= \ker(\chi)$.
\end{lemma}

\begin{proof}
Let $\chi\in\Irr(G)$ be afforded by a representation $\Phi$. From \cite[Lemma 2.27]{Is} we get that $\ZZ(\chi)=\{ g\in G : \Phi(g)= \varepsilon I \text{ for some }\varepsilon \in \C \}$.
Therefore, if there exists $g\in\ZZ(\chi)$ such that $\Phi(g)=\varepsilon I$ with $\varepsilon\ne 1$, that is $\ker(\chi) \ne \ZZ(\chi)$, then $\chi$ is not unisingular.
\end{proof}

\begin{prop}\label{centerless}
A  finite nonabelian group  $G$ such that $\ZZ(G)\neq 1$ always admits an irreducible nonlinear character which is not unisingular.
\end{prop}

\begin{proof}
Denote by $\Irr_1(G)$ the set of the irreducible nonlinear characters of $G$.
By Lemma \ref{ZK} it suffices to prove that there exists $\chi \in \Irr_1(G)$ such that $\ker(\chi) \neq \ZZ(\chi)$.

For the sake of contradiction, suppose that $\ker(\chi)= \ZZ(\chi)$ for any $\chi \in \Irr_1(G)$. By \cite[Theorem~$4.35$]{BZ}  we have
$$1=\bigcap_{\chi \in \Irr_1(G)} \ker(\chi)
=\bigcap_{\chi \in \Irr_1(G)} \ZZ(\chi)=\ZZ(G),$$
a contradiction.
\end{proof}

Note  that, if $\ZZ(G) \nsubseteq \ker(\chi)$, then $\ker(\chi)\ne\ZZ(\chi)$ and so the character $\chi$ is not unisingular.
Therefore, since $
\Irr\left(G/\ZZ(G)\right)=\left\{\chi\in\Irr(G):\ZZ(G)\subseteq\ker(\chi) \right\}$,
the irreducible unisingular characters of $G$ are exactly the irreducible unisingular characters of $G/\ZZ(G)$.
In conclusion, we can limit our study to the groups with trivial center.
In particular, since we are interested in representations of finite simple groups, in the following we will consider almost simple groups instead of quasisimple groups.

\section{The linear groups $\PSL_2(q)$}

In this section we study the unisingularity of the characters of the groups $\PSL_2(q)$ and $\PGL_2(q)$.
The main tool is Equation~\eqref{MolteplicityOf1}.
Let us start with the group $\PGL_2(q)$.
We must treat separately the case $q$ even and $q$ odd.

\subsection{The group $\PGL_2(q)$ with $q$ even}\label{PGL(2,q)even}

The character table of $\SL_2(q)$ is well known and for its construction we refer to \cite[Chapter 38]{Do}. As in that book, we set 
$$	I=\begin{pmatrix}
		1 & 0 \\
		0 & 1 \\
	\end{pmatrix},	
	\quad
	c=\begin{pmatrix}
		1 & 0 \\
		1 & 1 \\
	\end{pmatrix},
	\quad
    a=\begin{pmatrix}
		\nu & 0 \\
		0 & \nu^{-1} \\
	\end{pmatrix}$$
in $\SL_2(q)$, with $\nu$ being a generator of the multiplicative group of the finite field $\F_q$ of order $q$. Moreover, let $b\in \SL_2(q)$ be an element of order $q+1$.

Since $\PGL_2(q)\cong \SL_2(q)$, the character table of $\PGL_2(q)$ is described in Table~\ref{SL2even}, where 
$\rho \in \mathbb{C}$ is a primitive $(q-1)$-th root of $1$ and $\sigma \in \mathbb{C}$ is a primitive $( q+1)$-th root of~$1$, 
$1 \le i \le \frac{q-2}{2}$, $1 \le j \le \frac{q}{2}$, $1 \le \ell \le \frac{q-2}{2}$ and $1 \le m \le \frac{q}{2}$.
Note that $a$ is conjugate to $a^{-1}$ and $b$ is conjugate to $b^{-1}$, see \cite{Do}.

\begin{table}[ht]
	$\begin{array}{c|c c c c}
		& I & c & a^\ell & b^m  \\ \hline
		1_G & 1 & 1 & 1 & 1 \\\hline
		\psi & q &  0 & 1 & -1 \\\hline
		\chi_i & q+1 &  1 & \rho^{i\ell}+\rho^{-i\ell} & 0 \\\hline
		\theta_j & q-1 &  -1 & 0 & -\left(\sigma^{jm}+\sigma^{-jm}\right) \\
	\end{array} $
	\caption{Character table of $\PGL_2(q)\cong \SL_2(q)$ with $q$ even.}
	\label{SL2even}
\end{table}

In this case, $|c|=2$, $|a|=q-1$ and $|b| =q+1$. 
Since $a^\ell$ and $b^m$ are power of $a$ and $b$, to find the unisingularity of every character $\chi\in\Irr(\PGL_2(q))$ we will consider only $M_\chi(a)$ and $M_\chi(b)$.
Note that there are no other intersections between powers of elements belonging to different conjugacy classes.
By Equation~\eqref{MolteplicityOf1}, we have
$M_\chi(c)= \frac{1}{2}\left[\chi(I)+\chi(c)\right]$,
$M_\chi(a)= \frac{1}{q-1}\left[\chi(I)+\sum\limits_{t=1}^{q-2}\chi(a^t)\right]$
and $M_\chi(b)= \frac{1}{q+1}\left[\chi(I)+\sum\limits_{t=1}^{q}\chi(b^t)\right]$,
for any $\chi\in\Irr(\PGL_2(q))$.

We shall ignore the case where all the character values are nonnegative integers. Then, for the character $\psi$, we only have to consider
$M_\psi(b)=\frac{q-q}{q+1}=0$.
Hence, $\psi$ is not a unisingular character.

For the character $\chi_i$ we only have to compute
$$M_{\chi_i}(a)=\frac{1}{q-1}\left[q+1+\sum\limits_{t=1}^{q-2}\left( \rho^{it}+\rho^{-it} \right)\right].$$
Note that, for any $\varepsilon\neq 1$ that is an $f$-root of $1$, we have
\begin{equation}\label{n-rootSum}
	\sum_{n=1}^{f-1}\varepsilon^n = \frac{\varepsilon^f-1}{\varepsilon-1}-1=-1.
\end{equation}
Thus, $\sum\limits_{t=1}^{q-2}\left( \rho^{it}+\rho^{-it} \right)=-2$ and $M_{\chi_i}(a)=1$.
Therefore, $\chi_i$ is unisingular.

Lastly, for $\theta_j$ we have to compute $M_{\theta_j}$ for $c$ and $b$.
We obtain $M_{\theta_j}(c)=\frac{q-2}{2}$, which is equal to $0$ only if $q=2$.
For the element $b$,  we have
$$
M_{\theta_j}(b)=\frac{1}{q+1}\left[q-1-\sum\limits_{t=1}^{q}\left( \sigma^{jt}+\sigma^{-jt} \right)\right]
$$
and by \eqref{n-rootSum} we have $M_{\theta_j}(b)=1$.
Thus, the character $\theta_j$ is unisingular if and only if $q\ne 2$. 
Indeed, in this case it is a linear character.

\subsection{The group $\PGL_2(q)$ with $q$ odd}

The character table of $\PGL_2(q)$, when $q$ is a power of an odd prime $p$, can be deduced from that of $\GL_2(q)$, which can be found, for instance, in \cite{JL}.
In this case, the center $\ZZ$ of $\GL_2(q)$ consists of scalar matrices.

Let $a$, $b$, $c$ as before.
The character table of $\PGL_2(q)$ is described in Table~\ref{PGL2odd}, where 
\[
0\le k \le 1, \quad 1\le \ell \le \frac{q-3}{2}, \quad 1 \le m \le \frac{q-1}{2},\quad
\quad 1\le i \le \frac{q-3}{2} \equad 1\le j \le \frac{q-1}{2}.
\]
A complete set of representatives for the conjugacy classes consists of $\ZZ I$, $\ZZ c$, $\ZZ a^{\ell}$, $\ZZ a^{\frac{q-1}{2}}$, $\ZZ b^{m}$ and $\ZZ b^{\frac{q+1}{2}}$. Note that $|\ZZ c|=p$, $|\ZZ a|=q-1$, $|\ZZ b|=q+1$ and $|\ZZ a^{\frac{q-1}{2}}|=|\ZZ b^{\frac{q+1}{2}}|=2$.
Moreover, for each character the values of $\ZZ a^{\frac{q-1}{2}}$ and $\ZZ b^{\frac{q+1}{2}}$ are equal to the values of the families $\ZZ a^\ell$  with  $\ell = \frac{q-1}{2}$  and $\ZZ b^m$  with $m = \frac{q+1}{2}$. 
Hence, we can treat them together.

\begin{table}[htp]
	\[
	\begin{array}{c|cccccc}
		& \ZZ I & \ZZ c & \ZZ a^{\frac{q-1}{2}} & \ZZ b^{\frac{q+1}{2}} & \ZZ a^\ell & \ZZ b^m  \\ \hline
		\lambda_{k} & 1 & 1 & (-1)^{\frac{q-1}{2} k} &  (-1)^{\frac{q+1}{2}k} 
		& (-1)^{k\ell}& (-1)^{ m k}\\ \hline
		
		\psi_k & q   &  0   &  (-1)^{\frac{q-1}{2}k} & -(-1)^{\frac{q+1}{2} k} & (-1)^{k\ell}&
		-(-1)^{ m k}\\ \hline
		
		\chi_i & q+1 & 1 & 2(-1)^{ i} &0 & \rho^{i\ell}+\rho^{-i\ell} &0\\\hline 
		
		\theta_j & q-1 & -1 & 0 & -2(-1)^{ j} & 0 &
		-(\sigma^{ jm }+\sigma^{- jm})
	\end{array} 
	\]
	\caption{Character table of $\PGL_2(q)$, $q$ odd.}\label{PGL2odd}
\end{table}

Similarly to the case $q$ even, for any $\chi\in\Irr(\PGL_2(q))$ we have to compute $M_\chi$ only for the elements $\ZZ c$, $\ZZ a$ and $\ZZ b$, and their computations are the same as before.
We remark that, in this case $M_\chi(\ZZ c)=\frac{\chi(1)+(p-1)\chi(\ZZ c)}{p}$ for any $\chi \in \Irr(\PGL_2(q))$.
Note also that the multiplicity of the eigenvalue $1$ for the characters $\psi_0$, $\chi_i$ and $\theta_j$ are computed in the same way as  the characters $\psi$, $\chi_i$ and $\theta_j$ in the case $q$ even.
Note that in this case 

\[
M_{\theta_j}(\ZZ c)=\frac{q-1-(p-1)}{p}=\frac{q-p}{p}.
\]
Thus, $\chi_i$ is unisingular, $\psi_0$ is not unisingular and $\theta_j$ is unisingular if and only if $q$ is an odd prime.
We also observe that $\lambda_1$ is a non principal linear character and so it is not unisingular.

We only have left to compute $M_{\psi_1}$ on $\ZZ a$ and $\ZZ b$.
By Equation~\eqref{n-rootSum} we obtain 
$$
M_{\psi_1}(\ZZ a)= \frac{1}{q-1}
\left[q+\sum_{t=1}^{q-2}(-1)^{t}\right] =1 \equad 
M_{\psi_1}(\ZZ b)= \frac{1}{q+1}
\left[q-\sum_{t=1}^{q-2}(-1)^{t}\right]=1.
$$
Therefore $\psi_1$ is unisingular.

We can resume the results of these sections as follows (see also \cite[Theorem 1.1]{Z16}).
\begin{prop}\label{pPGL2}
	Let $\chi$ be an irreducible character of $\PGL_2(q)$. Then $\chi$ is not unisingular if and only if one of the following cases occurs:
	\begin{itemize}
		\item[(1)] $\chi$ is a nontrivial linear character;
        
		\item[(2)] $\chi(1) = q$ and $\chi(g)=1$, where $g$ has order $q-1$;
    
		\item[(3)] $q$ is an odd prime and $\chi(1)=q-1$.
	\end{itemize}
\end{prop}

Note that the characters described in the previous proposition at items $(2)$ and $(3)$ fail to be unisingular on elements of order $q+1$ and $p$, respectively.
Also, note that the character of item $(2)$ is the Steinberg character  of $\PGL_2(q)$.

\subsection{The group $\PSL_2(q)$}

Let now study $\PSL_2(q)$ when $q=p^k$ is a power of an odd prime $p$. 
We denote by $c$, $a$ and $b$ the matrices introduced in Section~\ref{PGL(2,q)even}.
Moreover, let us consider the following matrices in $\SL_2(q)$:
$$	z=\begin{pmatrix}
		-1 & 0 \\
		0 & -1 \\
	\end{pmatrix},
	\quad
	d=\begin{pmatrix}
		1 & 0 \\
		\nu & 1 \\
	\end{pmatrix}.$$
The center $\ZZ$ of $\SL_2(q)$ consists of $\{I,z\}$.
We can derive the character table of 
$\PSL_2(q)=\SL_2(q)/\ZZ$ from the character table of $\SL_2(q)$: we obtain two different cases
depending on the congruence class of $q$ modulo $4$.

When $q \equiv 1 \pmod{4}$ the character table is 
described in Table~\ref{PSL2oddCong1}, where
$$1 \le i \le \frac{q-5}{4},\quad 1 \le j \le \frac{q-1}{4}, \quad 1 \le \ell \le \frac{q-1}{4}\equad  1 \le m \le \frac{q-1}{4}.$$
When  $q \equiv 3\pmod{4}$ the character table is
described in Table~\ref{PSL2oddCong3}, where 
$$1 \le i \le \frac{q-3}{4},\quad 1 \le j \le \frac{q-3}{4},\quad 1 \le \ell \le \frac{q-3}{4}
\equad 1 \le m \le \frac{q+1}{4}.$$

\begin{table}[ht]
	\[
	\begin{array}{c|ccccc}
		& \ZZ I & \ZZ c & \ZZ d & \ZZ a^\ell & \ZZ b^m  \\ \hline
		1_G & 1 & 1 & 1 & 1 & 1 \\\hline
		\psi & q &  0 & 0 & 1 & -1 \\\hline
		\chi_i & q+1 &  1 & 1 & \rho^{2i\ell}+\rho^{-2i\ell} & 0 \\\hline
		\theta_j & q-1 &  -1 & -1 & 0 & -\left(\sigma^{2jm}+\sigma^{-2jm}\right) \\\hline
		\xi_1 & \frac{q+1}{2} & \frac{1+\sqrt{q}}{2} & \frac{1-\sqrt{q}}{2} & \left( -1\right)^\ell  & 0  \\\hline
		\xi_2 & \frac{q+1}{2} & \frac{1-\sqrt{q}}{2} & \frac{1+\sqrt{q}}{2} & \left( -1\right)^\ell  & 0  \\
	\end{array} 
	\]
	\caption{Character table of $\PSL_2(q)$ with $q \equiv 1 \pmod{4}$.}
	\label{PSL2oddCong1}
\end{table}

\begin{table}[htp]
	\[
	\begin{array}{c|ccccc}
		& \ZZ I & \ZZ c & \ZZ d & \ZZ a^\ell & \ZZ b^m  \\ \hline
		1_G & 1 & 1 & 1 & 1 & 1 \\\hline
		\psi & q &  0 & 0 & 1 & -1 \\\hline
		\chi_i & q+1 &  1 & 1 & \rho^{2i\ell}+\rho^{-2i\ell} & 0 \\\hline
		\theta_j & q-1 &  -1 & -1 & 0 & -\left(\sigma^{2jm}+\sigma^{-2jm}\right) \\\hline
		\eta_1 & \frac{q-1}{2} & \frac{-1+\sqrt{-q}}{2} & \frac{-1-\sqrt{-q}}{2} & 0 & \left( -1\right)^{m+1}  \\\hline
		\eta_2 & \frac{q-1}{2} & \frac{-1-\sqrt{-q}}{2} & \frac{-1+\sqrt{-q}}{2} & 0 & \left( -1\right)^{m+1}  \\
	\end{array} 
	\]
	\caption{Character table of $\PSL_2(q)$ with $q \equiv 3 \pmod{4}$.}\label{PSL2oddCong3}
\end{table}

A complete set of representatives for the conjugacy classes of $\PSL_2(q)$ consists of $\ZZ I$, $\ZZ c$, $\ZZ d$, $\ZZ a^\ell$ and $\ZZ b^m$. 
Note that $|\ZZ c|=|\ZZ d|=p$, $|\ZZ a|=\frac{q-1}{2}$ and 
$|\ZZ b|=\frac{q+1}{2}$. 
Given a positive integer $h$ such that $p \nmid h$, $(\ZZ c)^h$ belongs to the conjugacy class of $\ZZ c$ if and only if 
$h 1_{\F_q}$ is a square in $\F_q^\ast$. 
If this does not happen, then $(\ZZ c)^h$ belongs to the conjugacy class of $\ZZ d$, as $(\ZZ c)^h$ has order $p$.
Similarly, given a positive integer $t$ such that $p \nmid t$, $(\ZZ d)^t$ belongs to the conjugacy class of $\ZZ d$ if and only if 
$t 1_{\F_q}$ is a square in  $\F_q^\ast$.  
If this does not happen, then $(\ZZ d)^t$ belongs to the conjugacy class of $\ZZ c$. Hence, if $k$ is even, all the nontrivial powers of $\ZZ c$  belong to the same conjugacy class (the one of $\ZZ c$); otherwise, $\frac{p-1}{2}$ of these belong to the class of $\ZZ c$ and the other $\frac{p-1}{2}$ are conjugate to $\ZZ d$. The powers of $\ZZ d$  exhibit the same behavior. 
Therefore, for the computation we can consider only one of $M_{\chi}(\ZZ c)$ or $M_{\chi}(\ZZ d)$ when $k$ is odd, where $\chi \in \Irr(\PSL_2(q))$.

Since $|\PGL_2(q):\PSL_2(q)|=2$, for any $\chi \in \Irr(\PSL_2(q))$  the inertia group of $\chi$ is either $\PSL_2(q)$ or 
$\PGL_2(q)$.
By Clifford Theorem, any character of $\PGL_2(q)$ restricted to $\PSL_2(q)$ is either irreducible  or splits into $2$ irreducible characters of $\PSL_2(q)$, conjugated in $\PGL_2(q)$.
By degree reasons, $\psi$ is the restriction of both the characters of degree $q$ of $\PGL_2(q)$ and so it must be unisingular (since one of them is). 
Moreover, any $\theta_j$ is the restriction of a character of order $q-1$ with $j\ne \frac{q-1}{2}$. In $\PGL_2(q)$ they fail to be unisingular on the element $\ZZ c$ when $q$ is a prime. Thus, also in this case they are not unisingular when $q$ is a prime.
Hence, we only have to study the characters of degree $\frac{q+1}{2}$ when $q \equiv 1 \pmod4$ and $\frac{q-1}{2}$ when $q \equiv 3 \pmod4$.

Let us start by considering the case when $q \equiv 1 \pmod{4}$.
We have to compute $M_{\xi_1}$ only on $\ZZ d$ and $\ZZ a$. Note that in doing this we study also the unisingularity of $\xi_2$; indeed exchanging $\ZZ c$ and $\ZZ d$ we have the same computations.
We get

$$
M_{\xi_1}(\ZZ a)= \frac{2}{q-1}\left(\frac{q+1}{2}+\sum_{t=1}^{\frac{q-1}{2}-1}(-1)^{t}\right).$$
Since $\frac{q-1}{2}$ is even, $-1$ is a $\frac{q-1}{2}$-root of unity. Therefore, by Equation~\eqref{n-rootSum},
$$
M_{\xi_1}(\ZZ a)= \frac{2}{q-1}\left(\frac{q+1}{2}-1\right) =1.
$$
For the conjugacy class $\ZZ d$ we need to distinguish two cases, according to the parity of $k$. If $k$ is odd, then 
\[
	M_{\xi_1}(\ZZ d) = \dfrac{1}{p}\left( \dfrac{q+1}{2}+\dfrac{1+\sqrt{q}}{2}\cdot \dfrac{p-1}{2} + \dfrac{1-\sqrt{q}}{2} \cdot\dfrac{p-1}{2} \right) =\dfrac{p^{k-1}+1}{2}\geq 1.
\]
If $k$ is even, then
\[
M_{\xi_1}(\ZZ d)= \dfrac{1}{p}\left( \dfrac{q+1}{2}+\left(p-1\right) \dfrac{1-\sqrt{q}}{2} \right)=\frac{p^{k-1}-p^\frac{k}{2}+p^{\frac{k}{2}-1} +1}{2}\geq 1.
\]

Finally, let us consider $q \equiv 3 \pmod{4}$. 
We must compute $M_{\eta_1}$ for $\ZZ c$ and $\ZZ b^m$. Indeed, since $q \equiv 3 \pmod{4}$, $k$ must be odd and we do not have to consider $\ZZ d$ as it is a power of $\ZZ c$.
As for $\xi_2$, we do not need to study $\eta_2$.
We obtain 
\[M_{\eta_1}(\ZZ b)=\frac{2}{q+1}\left(\frac{q-1}{2}-\sum_{t=1}^{\frac{q+1}{2}-1}(-1)^{t}\right)= 1\]
and  
\[M_{\eta_1}(\ZZ c)=\dfrac{1}{p}\left( \dfrac{q-1}{2}+\dfrac{-1+\sqrt{-q}}{2}\cdot \dfrac{p-1}{2} + \dfrac{-1-\sqrt{-q}}{2} \cdot\dfrac{p-1}{2} \right) = \frac{p^{k-1}-1}{2}.\]
The last value is positive if and only if $q$ is not a prime.

We can resume the results of this section as follows (see also 
\cite[Theorem 5.3]{Z16} and \cite[Proposition 10]{Z25}). 
\begin{prop}
	Let $\chi$ be an irreducible character of $\PSL_2(q)$.
	Then $\chi$ is not unisingular if and only if one of the following cases occurs:
	\begin{itemize}
		\item[(1)] $q\in \{2,3\}$ and $\chi$ is a nontrivial linear character;
		\item[(2)] $q$ is even and $\chi(1)=q$;
		\item[(3)] $q$ is an odd prime and $\chi(1)=q-1$;
		\item[(4)] $q>3$ is an odd prime such that $q\equiv 3 \pmod 4$ and $\chi(1)=\frac{q-1}{2}$.
	\end{itemize}
\end{prop}

Note that the characters at items (2), (3) and (4) fail to be unisingular on elements of order $q+1$, $p$ and $p$, respectively. 

\section{The unitary groups $\PSU_3(q)$}\label{U3}

In this section we study the unisingularity of the characters of 
$\PGU_3(q)$ and $\PSU_3(q)$, where $q$ is a power of a prime $p$.
The character table of $\GU_3(q)$ was described in \cite{En}.
Let $\eta$ be an element of order $q^2-1$ in $\F_{q^2}^*$ and $\lambda$ be an element of order $q^3+1$ in $\F_{q^6}^*$; take an element $\nu$ of order $q+1$ in $\F_{q^2}^*$ such that $\nu=\eta^{q-1}=\lambda^{q^2-q+1}$.
To describe the conjugacy classes of $\GU_3(q)$, 
let us consider the following matrices of $\GL_3(q^6)$:
$$\overline{A} = \begin{pmatrix}
		1 & 0 & 0 \\
		1 & 1 & 0 \\
		0 & 0 & 1 \\
	\end{pmatrix},\quad
    \overline{B} = \begin{pmatrix}
		1 & 0 & 0 \\
		1 & 1 & 0 \\
		0 & 1 & 1 \\
	\end{pmatrix},\quad
        \overline{C} = \begin{pmatrix}
		\nu & 0 & 0 \\
		0 & 1 & 0 \\
		0 & 0 & 1 \\
	\end{pmatrix},\quad
    \overline{D_b} = \begin{pmatrix}
		\nu^b & 0 & 0 \\
		0 & 1 & 0 \\
		0 & 1 & 1 \\
	\end{pmatrix},$$
    $$
        \overline{E_{\alpha,\beta}} = \begin{pmatrix}
		\nu^\alpha & 0 & 0 \\
		0 & \nu^\beta & 0 \\
		0 & 0 & 1 \\
	\end{pmatrix},
\quad
    \overline{F}= \begin{pmatrix}
        \eta & 0 & 0\\
        0 & \eta^{-q} & 0 \\
        0 & 0 & 1
    \end{pmatrix},
\quad
    \overline{G}=\begin{pmatrix}
    \lambda & 0 & 0\\
    0 & \lambda^{-q} & 0\\
    0 & 0 & \lambda^{q^2}
\end{pmatrix},$$
where $1\le b\le q$ and $1 \le \alpha < \beta \le q$.
Let $A,B,C,D_b,E_{\alpha,\beta},F,G$ be elements of $\GU_3(q)$ conjugated respectively to the matrices $\overline{A},\overline{B}, \overline{C},\overline{D_b},\overline{E_{\alpha,\beta}}, \overline{F},\overline{G}$ in $\GL_3(q^6)$. The center $\ZZ$ of $\GU_3(q)$ consists of the $q+1$  scalar matrices 
$\nu^r I$, where $0\le r\le q$.

Let $\mathscr{C}$ be a complete set of representatives for the quotient set $\{c : 1\le c\le q^2-1, (q-1)\nmid c\}/\sim$, where
\[
c\sim \overline{c} \Leftrightarrow \overline{c}\equiv c,-cq \pmod{q^2-1}.
\]
Let $\mathscr{D}$ be a complete set of representatives for the quotient set $\{d: 1\le d\le q^2-q\}/\vartriangle$, where
\[
d \vartriangle \overline{d} \Leftrightarrow \overline{d}\equiv d,-dq,dq^2 \pmod{q^2-q+1}.
\]
Let $\mathscr{E}$ be a complete set of representatives for the quotient set $\{(\alpha,\beta):1\le\alpha<\beta\le q\}/ \diamond$, where
\[
(\alpha,\beta)\diamond (\gamma,\delta) \Leftrightarrow \{\nu^\gamma,\nu^\delta\}\in\{ \{\nu^\alpha,\nu^\beta\}, \{\nu^{-\alpha},\nu^{\beta-\alpha}\}, \{\nu^{\alpha-\beta},\nu^{-\beta}\} \}.
\]

\subsection{The group $\PGU_3(q)$ with $\gcd(q+1,3)=1$}

First of all, note that $|\mathscr{C}|=\frac{(q+1)(q-2)}{2}$, $|\mathscr{D}|=\frac{q(q-1)}{3}$ and $|\mathscr{E}|=\frac{q(q-1)}{6}$.

A complete set of representatives for $\PGU_3(q)$, when $3\nmid(q+1)$, consists of 
\[
\ZZ I,\quad \ZZ A,\quad \ZZ B,\quad (\ZZ C)^a,\quad \ZZ D_b,\quad \ZZ E_{\alpha,\beta},\quad \ZZ F^c \equad \ZZ G^d
\]
where $
a,b\in\{1,\dots,q\}$, $c\in \mathscr{C}$, $d\in\mathscr{D}$ and $(\alpha,\beta)\in \mathscr{E}$.
The character table is given in Table~\ref{PGU3n1}, where $\rho$ is a primitive $(q+1)$-root of $1$, $\sigma$ is a primitive $(q^2-1)$-root of $1$, $\tau$ is a primitive $(q^2-q+1)$-root of $1$, $1\le h\le q$, $i \in \mathscr{C}$, $j\in \mathscr{D}$ and $(m,n)\in \mathscr{E}$.

\begin{table}[tph]
	\rotatebox{90}{
		\begin{footnotesize}
			$\begin{array}{c|ccccc}
				
				& \ZZ I & \ZZ A & \ZZ B & (\ZZ C)^a & \ZZ D_b  \\ \hline

				\chi_{1}   & 1 & 1 & 1 & 1 & 1 \\\hline
				
				\chi_{q^2-q}   & q^2-q  & -q & 0 & -(q-1) & 1 \\\hline
				
				\chi_{q^3}   & q^3 & 0 & 0 & q & 0 \\\hline
				
				\chi_{q^2-q+1}^{(h)}   & q^2-q+1 & -(q-1) & 1 & -(q-1)\rho^{ah}+\rho^{-2ah} & \rho^{bh}+\rho^{-2bh} \\\hline
				
				\chi_{q(q^2-q+1)}^{(h)}   & q(q^2-q+1) & q & 0 & (q-1)\rho^{ah}+q\rho^{-2ah} & -\rho^{bh} \\\hline
				
				\chi_{(q-1)(q^2-q+1)}^{(n,m)}   & 
                      (q-1)(q^2-q+1) & 2q-1 & -1 & 
                    \begin{array}{c}(q-1)\rho^{a(n-2m)}+\\
                      (q-1)(\rho^{a(m-2n)}+\rho^{a(n+m)})
                      \end{array}
                      & 
                                          -\rho^{b(n-2m)}-\rho^{b(m-2n)}
                      -\rho^{b(n+m)}\\\hline
				
				\chi_{(q+1)(q^2-q+1)}^{(i)}   & (q+1)(q^2-q+1) & 1 & 1 & (q+1)\rho^{ai} & \rho^{bi} \\\hline
				
				\chi_{(q-1)(q+1)^2}^{(j)}   & (q-1)(q+1)^2 & -(q+1) & -1 & 0 & 0 \\
				
			\end{array}$
	\end{footnotesize}        }
    \quad
	\rotatebox{90}{
		\begin{footnotesize}
			$\begin{array}{c|cccc}
				
				& \ZZ E_{\alpha,\beta} & (\ZZ F)^c & (\ZZ G)^d  \\ \hline

				\chi_{1}   & 1 & 1 & 1\\\hline
				
				\chi_{q^2-q}   & 2 & 0 & -1 \\\hline
				
				\chi_{q^3}   & -1 & 1 & -1 \\\hline
				
				\chi_{q^2-q+1}^{(h)}   & \rho^{\alpha h}+ \rho^{\beta h} + \rho^{-(\alpha+\beta )h} & \rho^{ch} & 0 \\\hline
				
				\chi_{q(q^2-q+1)}^{(h)}  & -\rho^{\alpha h}- \rho^{\beta h} - \rho^{-(\alpha +\beta)h} & \rho^{ch} & 0 \\\hline
				
				\chi_{(q-1)(q^2-q+1)}^{(n,m)}   &  
                \begin{array}{c}
                -\rho^{\alpha(n-m)-\beta m}-\rho^{\alpha(m-n)-\beta n}-\rho^{\beta(n-m)-\alpha m}-\rho^{\beta(m-n)-\alpha n}\\
                - \rho^{\alpha n+\beta m}-\rho^{\alpha m+\beta n} 
                \end{array}
                & 0 & 0 \\\hline
				
				\chi_{(q+1)(q^2-q+1)}^{(i)}   & 0 & \sigma^{ci}+\sigma^{-qci} & 0 \\\hline
				
				\chi_{(q-1)(q+1)^2}^{(j)}   & 0 & 0 & -\tau^{dj}-\tau^{-qdj}-\tau^{q^2dj} \\
				
			\end{array}$
	\end{footnotesize}}
	
	\caption{Character table of $\PGU_3(q)$, with $\gcd(q+1,3)=1$.}
	\label{PGU3n1}
	
\end{table}

For the classes $\ZZ E_{\alpha,\beta}$ and $\ZZ D_b$ we can limit to study $\ZZ E_{\alpha,\beta}$ when $\gcd(\alpha,\beta)=1$ and $\ZZ D_1$. 
Moreover, we obtain that the orders of the elements in each conjugacy class are
$$o(\ZZ I) = 1, \quad o(\ZZ A) = p, \quad o(\ZZ B) = \begin{cases}
	4 & \text{if } p=2 \\
	p & \text{otherwise}
\end{cases},\quad
o(\ZZ C)= q+1,$$
$$o(\ZZ D_1)= p(q+1), \quad o(\ZZ E_{\alpha,\beta})= q+1,\quad 
o(\ZZ F)= q^2-1, \quad o(\ZZ G) =	q^2-q+1,$$
and for the power maps we have: 
\begin{itemize}
    \item the element $(\ZZ A)^k$ is conjugate to $\ZZ A$ for any $1\le k\le p-1$; 
    \item if $p\ne 2$ then $(\ZZ B)^k$ is conjugate to $\ZZ B$ for any $1\le k\le p-1$; if $p=2$, then $(\ZZ B)^3$ is conjugate to $\ZZ B$, while $(\ZZ B)^2$ is conjugate to $\ZZ A$;
    \item if $p \mid k$, then $(\ZZ D_1)^k$ is conjugate to $(\ZZ C)^k$; if $(q+1) \mid k$, then $(\ZZ D_1)^k$ is conjugate to $\ZZ A$;
    \item 
    if $k\alpha  \equiv0 \pmod{q+1}$, then $(\ZZ E_{\alpha,\beta})^k$ is conjugate to $(\ZZ C)^{k\beta}$;
    if $k\beta \equiv0 \pmod{q+1}$, then $(\ZZ E_{\alpha,\beta})^k$ is conjugate to $(\ZZ C)^{k\alpha}$;
    if $k\alpha\equiv k\beta \pmod{q+1}$, then  $(\ZZ E_{\alpha,\beta})^k$ is conjugate to $(\ZZ C)^{-k\alpha}$;
    \item if $(q-1) \mid c$, then $(\ZZ F)^c$ is  conjugate to $(\ZZ C)^{-\frac{c}{q-1}}$.
\end{itemize}
Furthermore, there are no other intersections between powers of elements belonging to different conjugacy classes.
So, we can limit our study to the elements $\ZZ B$, $\ZZ D_1$, $\ZZ E_{\alpha,\beta}$, $\ZZ F$ and $\ZZ G$.

Now, let $g_\alpha=\gcd(\alpha,q+1)$, $g_\beta=\gcd(\beta,q+1)$ and $\chi\in\Irr(\PGU_3(q))$ with $3 \nmid (q+1)$. 
Also, let $\gamma=\beta-\alpha$ and $g_\gamma=\gcd(\gamma,q+1)$. We get
$$
M_\chi(\ZZ B)=\begin{cases}
\frac{1}{4}(\chi(\ZZ I)+ \chi(\ZZ A)+ 2 \chi(\ZZ B)) & \text{if } p=2, \\
\frac{1}{p}(\chi(\ZZ I)+(p-1)\chi(\ZZ B)) & \text{if } p\geq 3,    
\end{cases}$$
$$M_\chi(\ZZ D_1)  =  \frac{1}{p(q+1)}\left(\chi(\ZZ I)+ (p-1)\chi(\ZZ A)+ \sum\limits_{a=1}^{q}\chi(\ZZ C^a) +\sum\limits_{\substack{t=1 \\ p, (q+1)\nmid t}}^{p(q+1)}\chi(\ZZ D_t)\right),
$$
\[
\begin{array}{rcl}
M_\chi(\ZZ E_{\alpha,\beta}) & = &  \frac{1}{q+1}\left(\chi(\ZZ I)+
\sum\limits_{t=1}^{g_\alpha-1}\chi(\ZZ C^{t\frac{q+1}{g_\alpha}}) +\sum\limits_{t=1}^{g_\beta-1}\chi(\ZZ C^{t\frac{q+1}{g_\beta}}) \right)\\
& +  & \frac{1}{q+1}\left(\sum\limits_{t=1}^{g_\gamma-1}\chi(\ZZ C^{t\frac{q+1}{g_\gamma}})+ 
\sum\limits_{\substack{t=1 \\ \frac{q+1}{g_\alpha}, \frac{q+1}{g_\beta}, \frac{q+1}{g_\gamma}\nmid t}}^{q}\chi(\ZZ E_{t\alpha,t \beta })\right),
\end{array}
\]
\[
M_\chi(\ZZ F)=\frac{1}{q^2-1}\left(\chi(\ZZ I)+\sum_{a=1}^{q}\chi(\ZZ C^a) + \sum_{\substack{t=1 \\ (q-1) \nmid t}}^{q^2-1}\chi(\ZZ F^t)\right) ,
\]
\[
M_\chi(\ZZ G)=\frac{1}{q^2-q+1}\left(\chi(\ZZ I)+\sum_{t=1}^{q^2-q}\chi(\ZZ G^t) \right).
\]

Skipping the trivial cases, we start with the character $\chi_{q^2-q}$.
We can ignore $M_{\chi_{q^2-q}}(\ZZ B)$ when $p\ne2$.
When $p=2$, we have
\[
M_{\chi_{q^2-q}}(\ZZ B)=\frac{1}{4}\left[q^2-q-q\right]=\frac{q^2-2q}{4},
\]
that is positive if and only if $q> 2$.
For $\ZZ D_1$ we get:
$$
M_{\chi_{q^2-q}}(\ZZ D_1) = \frac{1}{p(q+1)}\left[q^2-q-(p-1)q-q(q-1)+p(q+1)-q-p \right]=0.$$
Continuing with $\ZZ E_{\alpha,\beta}$, we obtain

\[
\begin{array}{rcl}
M_{\chi_{q^2-q}}(\ZZ E_{\alpha,\beta})&=&\frac{1}{q+1}\left[q^2-q-(q-1)(g_\alpha+g_\beta+g_\gamma-3)\right. \\
&&\left.+2q-2(g_\alpha+g_\beta+g_\gamma-3) \right]\\
&=&\frac{1}{q+1}\left[q^2+q-(q+1)(g_\alpha+g_\beta+g_\gamma-3)\right].
\end{array}
\]
We remark that if $x,y$ are two nonnegative integers, then $x+y\le xy+1$ with equality holding if and only if $x$ or $y$ are equal to $1$. 
Indeed, we have $xy+1-x-y=(x-1)(y-1)\ge0$.
From this equality it is immediate to prove that if $x,y,z$ are three
nonnegative integer, then $x+y+z\le xyz+2$ with equality holding if and only if at least two of $x,y,z$ are equal to $1$.
Since $1\leq \alpha<\beta\leq q$, $1\le g_\alpha, g_\beta, g_\gamma\leq q$, 
$\gcd(g_\alpha,g_\beta) =\gcd(g_\alpha,g_\gamma)=\gcd(g_\beta,g_\gamma)=1$ and $g_\alpha g_\beta g_\gamma\mid q+1$, we have $g_\alpha+g_\beta+g_\gamma<q+3$, whence
$$
M_{\chi_{q^2-q}}(\ZZ E_{\alpha,\beta})>\frac{q^2+q-(q+1)q}{q+1}=0.$$
For the remaining classes, we find
$$
M_{\chi_{q^2-q}}(\ZZ F)=\frac{q^2-q-q(q-1)}{q^2-1}=0
\equad 
M_{\chi_{q^2-q}}(\ZZ G)=\frac{q^2-q-(q^2-q)}{q^2-q+1}=0.$$
We conclude that the character $\chi_{q^2-q}$ is not unisingular.

For the character $\chi_{q^3}$ we can limit to compute $M_{\chi_{q^3}}$ for $\ZZ E_{\alpha,\beta}$ and $\ZZ G$.
We get
\[
\begin{array}{rcl}
M_{\chi_{q^3}}(\ZZ E_{\alpha,\beta})&=&\frac{1}{q+1}\left[q^3+q(g_\alpha+g_\beta+g_\gamma-3)-q-(g_\alpha+g_\beta+g_\gamma-3) \right]\\
&=&\frac{1}{q+1}\left[q^3-q+(q-1)(g_\alpha+g_\beta+g_\gamma-3)\right].
\end{array}
\]
Since $g_\alpha,g_\beta,g_\gamma\ge1$, then
\[
M_{\chi_{q^3}}(\ZZ E_{\alpha,\beta})\ge \frac{q^3-q}{q+1}=q^2-q>0.
\]
For $\ZZ G$, we obtain
\[
M_{\chi_{q^3}}(\ZZ G)=\frac{q^3-q^2+q}{q^2-q+1}=q>0.
\]
Hence, the character $\chi_{q^3}$ is unisingular.

For the character $\chi_{q^2-q+1}^{(h)}$ we can ignore the element $\ZZ G$. 
Moreover, it suffices to study $M_{\chi_{q^2-q+1}^{(h)}}(\ZZ B)$ only for $p=2$. 
In this case, we have
\[M_{
\chi_{q^2-q+1}^{(h)}}(\ZZ B)=\frac{1}{4}\left[q^2-q+1-(q-1)+2\right]=\frac{q^2-2q+4}{4}>0.
\]
For $\ZZ D_1$, called $M=M_{\chi_{q^2-q+1}^{(h)}}(\ZZ D_1)$ we get
\[
\begin{array}{rcl}
M&=&\frac{1}{p(q+1)}\left[q^2-q+1-(p-1)(q-1)
+\sum\limits_{t=1}^{q}(-(q-1)\rho^{ht}+\rho^{-2ht})\right.\\
&&\left.+\sum\limits_{t=1}^{p(q+1)}(\rho^{ht}+\rho^{-2ht})
-\sum\limits_{t=1}^{q}(\rho^{hpt}+\rho^{-2hpt})
-\sum\limits_{t=1}^{p}(\rho^{h(q+1)t}+\rho^{-2h(q+1)t})\right].
\end{array}
\]
Since $\rho^h$ and $\rho^{hp}$ are two primitive $(q+1)$-roots of $1$, we have

\[
\sum\limits_{t=1}^{q}(\rho^{ht})=\sum\limits_{t=1}^{q}(\rho^{hpt})=-1 \equad \sum\limits_{t=1}^{q}(\rho^{-2ht})=\sum\limits_{t=1}^{q}(\rho^{-2hpt}).
\]
Therefore,
\[
\begin{array}{rcl}
M&=&\frac{1}{p(q+1)}\left[q^2-pq+p+(q-1)
+\sum\limits_{t=1}^{q}(\rho^{-2ht})+\sum\limits_{t=1}^{p(q+1)}(\rho^{-2ht})
+1-\sum\limits_{t=1}^{q}(\rho^{-2ht})-2p\right]\\
&=& \frac{1}{p(q+1)}\left[q^2-pq-p+q+\sum\limits_{t=1}^{p(q+1)}(\rho^{-2ht})\right].
\end{array}
\]
Now, if  $q$ is even, then $2\nmid (q+1)$ and $\rho^{-2h}$ is a primitive $(q+1)$-root of $1$; so $q>2$ and
\[
M\ge\frac{1}{p(q+1)}\left[q^2-pq-p+q\right]=\frac{(q-p)(q+1)}{p(q+1)}=\frac{q-2}{2}>0.
\]
If $q$ is odd, then $2\mid (q+1)$ and so: if $h\ne \frac{q+1}{2}$ we have $\rho^{-2h}\ne 1$, whence $M=\frac{q-p}{p}$; if $h=\frac{q+1}{2}$ we have $\rho^{-2h}= 1$, whence
\[
M\ge\frac{1}{p(q+1)}\left[q^2-pq-p+q+p(q+1)\right]=\frac{(q-p)(q+1)}{p(q+1)}+1=\frac{q}{p}>0.
\]
We conclude that $M_{\chi_{q^2-q+1}^{(h)}}(\ZZ D_1)=0$ if and only if $q$ is an odd prime and $h\ne \frac{q+1}{2}$.
Continuing with $\ZZ E_{\alpha,\beta}$, let $\theta_\alpha$ be a primitive $g_\alpha$-root of $1$, $\theta_\beta$ be a primitive $g_\beta$-root of $1$ and $\theta_\gamma$ be a primitive $g_\gamma$-root of $1$.
Calling $M=M_{\chi_{q^2-q+1}^{(h)}}(E_{\alpha,\beta})$, we obtain
\[
\begin{array}{rcl}
M&=&\frac{1}{q+1}\left[q^2-q+1+\sum\limits_{t=1}^{g_\alpha-1}(-(q-1)\rho^{\frac{q+1}{g_\alpha}ht}+\rho^{-2\frac{q+1}{g_\alpha}ht})\right.\\
&&+\sum\limits_{t=1}^{g_\beta-1}(-(q-1)\rho^{\frac{q+1}{g_\beta}ht}+\rho^{-2\frac{q+1}{g_\beta}ht}) +\sum\limits_{t=1}^{g_\gamma-1}(-(q-1)\rho^{\frac{q+1}{g_\gamma}ht}+\rho^{-2\frac{q+1}{g_\gamma}ht})\\
&&+\sum\limits_{t=1}^{q}(\rho^{\alpha ht}+\rho^{\beta ht}+\rho^{-(\alpha+\beta) ht}-\sum\limits_{t=1}^{g_\alpha-1}(\rho^{\alpha\frac{q+1}{g_\alpha}ht}+\rho^{\beta\frac{q+1}{g_\alpha}ht}+\rho^{-(\alpha+\beta)\frac{q+1}{g_\alpha}ht})\\
&&-\sum\limits_{t=1}^{g_\beta-1}(\rho^{\alpha\frac{q+1}{g_\beta}ht}+\rho^{\beta\frac{q+1}{g_\beta}ht}+\rho^{-(\alpha+\beta)\frac{q+1}{g_\beta}ht}) \\
&&\left.  -\sum\limits_{t=1}^{g_\gamma-1}(\rho^{\alpha\frac{q+1}{g_\gamma}ht}+\rho^{\beta\frac{q+1}{g_\gamma}ht}+\rho^{-(\alpha+\beta)\frac{q+1}{g_\gamma}ht})\right]\\
&=&\frac{1}{q+1}\left[q^2-q+1+\sum\limits_{t=1}^{g_\alpha-1}(-(q-1)\theta_\alpha^{ht}+\theta_\alpha^{-2ht}) +\sum\limits_{t=1}^{g_\beta-1}(-(q-1)\theta_\beta^{ht}+\theta_\beta^{-2ht}) \right.\\
&& +\sum\limits_{t=1}^{g_\gamma-1}(-(q-1)\theta_\gamma^{ht}+\theta_\gamma^{-2ht}) +\sum\limits_{t=1}^{q}(\rho^{\alpha ht}+\rho^{\beta ht}+\rho^{-(\alpha+\beta) ht})\\
&& -\sum\limits_{t=1}^{g_\alpha-1}(\theta_\alpha^{\alpha ht}+\theta_\alpha^{\beta ht}+\theta_\alpha^{-(\alpha+\beta) ht})-\sum\limits_{t=1}^{g_\beta-1}(\theta_\beta^{\alpha ht}+\theta_\beta^{\beta ht}+\theta_\beta^{-(\alpha+\beta) ht})\\
&& \left.-\sum\limits_{t=1}^{g_\gamma-1}(\theta_\gamma^{\alpha ht}+\theta_\gamma^{\beta ht}+\theta_\gamma^{-(\alpha+\beta) ht})\right].
\end{array}
\]
Since $\theta_\alpha^\alpha=1$ and $\gcd(\beta,g_\alpha)=1$, we have that $\theta_\alpha^\beta$ is a primitive $g_\alpha$-root of $1$.
Also, $\theta_\gamma^\alpha$ is a primitive $g_\gamma$-root of $1$ and $\theta_\gamma^{\alpha+\beta}=\theta_\gamma^{2\alpha}$.
In particular,  $\sum\limits_{t=1}^{g_\alpha-1} \theta_\alpha^{-(\alpha+\beta) ht} = \sum\limits_{t=1}^{g_\alpha-1} \theta_\alpha^{ht}.$
Similarly, $\theta_\beta^\alpha$ is a primitive $g_\beta$-root of $1$ and $\sum\limits_{t=1}^{g_\beta-1} \theta_\beta^{-(\alpha+\beta) ht} = \sum\limits_{t=1}^{g_\beta-1} \theta_\beta^{ht}.$
Thus,
\[
\begin{array}{rcl}
M&=&\frac{1}{q+1}\left[q^2-q+3-g_\alpha-g_\beta+\sum\limits_{t=1}^{g_\alpha-1}(-(q-1)\theta_\alpha^{ht}+\theta_\alpha^{-2ht}) \right.\\
&&+\sum\limits_{t=1}^{g_\beta-1}(-(q-1)\theta_\beta^{ht}+\theta_\beta^{-2ht}) +\sum\limits_{t=1}^{g_\gamma-1}(-(q-1)\theta_\gamma^{ht}+\theta_\gamma^{-2ht}) \\
&&\left. +\sum\limits_{t=1}^{q}(\rho^{\alpha ht}+\rho^{\beta ht}+\rho^{-(\alpha+\beta) ht}) -2\sum\limits_{t=1}^{g_\alpha-1}(\theta_\alpha^{ht}) -2\sum\limits_{t=1}^{g_\beta-1}(\theta_\beta^{ht}) -\sum\limits_{t=1}^{g_\gamma-1}(2\theta_\gamma^{ ht}+\theta_\gamma^{-2 ht})\right]\\
&= &\frac{1}{q+1}\left[q^2-q+3-g_\alpha-g_\beta+\sum\limits_{t=1}^{g_\alpha-1}(-(q+1)\theta_\alpha^{ht}+\theta_\alpha^{-2ht})\right. \\
&&\left. +\sum\limits_{t=1}^{g_\beta-1}(-(q+1)\theta_\beta^{ht}+\theta_\beta^{-2ht}) -(q+1)\sum\limits_{t=1}^{g_\gamma-1}(\theta_\gamma^{ht}) +\sum\limits_{t=1}^{q}(\rho^{\alpha ht}+\rho^{\beta ht}+\rho^{-(\alpha+\beta) ht})\right].\\
\end{array}
\]
Note that $\sum\limits_{t=1}^{g_\alpha-1}(-(q+1)\theta_\alpha^{ht}+\theta_\alpha^{-2ht})$ assumes the value $q$ if $\theta_\alpha^{-2h}\ne 1$, the value $q+g_\alpha$ if $\theta_\alpha^{-2h}=1$ and $\theta_{\alpha}^h\ne1$ and the value $-(q+1)(g_\alpha-1)+(g_\alpha-1)$ if $\theta_\alpha^{-2h}=\theta_{\alpha}^h=1$.
The minimum value is $-(q+1)(g_\alpha-1)+(g_\alpha-1)$.
In the same way, the mimimum value of $\sum\limits_{t=1}^{g_\beta-1}(-(q+1)\theta_\beta^{ht}+\theta_\beta^{-2ht})$ is $-(q+1)(g_\beta-1)+(g_\beta-1)$.
Then
\[
\begin{array}{rcl}
M&\ge&\frac{1}{q+1}[q^2-q+3-g_\alpha-g_\beta-(q+1)(g_\alpha+g_\beta+g_\gamma-3)+g_\alpha+g_\beta-5]\\
&=&\frac{1}{q+1}[q^2-q-2-(q+1)(g_\alpha+g_\beta+g_\gamma-3)]. 
\end{array}
\]
We now need an auxiliary result.

\begin{lemma}\label{DisugualianzaSommaProdotto}
    Let $x,y,z,m$ be four positive integers such that $x,y,z\ne m$, $x\le y\le z$, $x,y,z$ are pairwise coprime and $xyz\mid m$.
    If $m>6$, then
    \[
    x+y+z<m.
    \]
\end{lemma}
\begin{proof}
    Let $m=xyzw$.
    Since $x+y+z\le3z$, if $xyw>3$ we have $x+y+z\le 3z<xyzw=m$.
    If $x=y=w=1$, then $z=m$ and so we ignore this case.
    Let us now study the remaining cases:
    \begin{enumerate}[$(1)$]
        \item if $(x,y,w)=(1,2,1)$, then  $x+y+z=3+z$ and $m=2z$.
        In this case, $x+y+z<m$ if and only if $z>3$, that is $m>6$;
        \item if $(x,y,w)=(1,1,2)$, then $x+y+z=2+z$ and $m=2z$.
        In this case,  $x+y+z<m$ if and only if $z>2$, that is $m>4$;
        \item if $(x,y,w)=(1,3,1)$, then $x+y+z=4+z$ and $m=3z$.
        In this case, $x+y+z<m$ if and only if $z>2$, that is $m>6$;
        \item if $(x,y,w)=(1,1,3)$, then $x+y+z=2+z$ and $m=3z$.
        In this case,  $x+y+z<m$ if and only if $z>2$, that is $m>6$. \qedhere
    \end{enumerate}
\end{proof}

Back to our problem, since $g_\alpha$, $g_\beta$, $g_\gamma$ and $q+1$ satisfy the hypotheses of Lemma~\ref{DisugualianzaSommaProdotto}, when $q>5$ we obtain
\[
M>\frac{1}{q+1}\left[q^2-q-2-(q+1)(q-2)\right]=0.
\]
We now consider the exceptional cases $q=3$ and $q=4$.
If $q=3$, then $|\mathscr{E}|=1$, $\ZZ E_{1,3}$ has order $4$ and
$M_{\chi_{q^2-q+1}^{(h)}}(\ZZ E_{1,3})\in \{1,3\}$.
If $q=4$, then  $|\mathscr{E}|=2$, $\ZZ E_{1,4}$ has order $5$ and
$M_{\chi_{q^2-q+1}^{(h)}}(\ZZ E_{1,4})=3$.
Lastly, for the element $\ZZ F$ we have
\[
\begin{array}{rcl}
M_{\chi_{q^2-q+1}^{(h)}}(\ZZ F)&=& \frac{1}{q^2-1}\left[q^2-q+1+\sum\limits_{t=1}^{q}(-(q-1)\rho^{ht}+\rho^{-2ht})\right.\\
&&\left.+\sum\limits_{t=1}^{q^2-2}(\rho^{th})-\sum\limits_{t=1}^{q}(\rho^{(q-1)th})\right]\\
&=&\frac{1}{q^2-1}\left[q^2-q+1+\sum\limits_{t=1}^{q}(-(q-1)\rho^{ht}+\rho^{-2ht})\right.\\
&&\left.+\sum\limits_{t=1}^{q^2-2}(\rho^{ht})-\sum\limits_{t=1}^{q}(\rho^{-2ht})\right]=1.
\end{array}
\]
In conclusion, the character $\chi_{q^2-q+1}^{(h)}$ is not unisingular if and only if $h\ne \frac{q+1}{2}$ and $q=p>2$.

For the character $\chi_{q(q^2-q+1)}^{(h)}$ we can compute the value of $M_{\chi_{q(q^2-q+1)}^{(h)}}$ only for the elements $\ZZ D_{1}$, $\ZZ E_{\alpha,\beta}$ and $\ZZ F$.
Calling $M=M_{\chi_{q(q^2-q+1)}^{(h)}}(\ZZ D_1)$, we get
\[
\begin{array}{rcl}
M&=&\frac{1}{p(q+1)}\left[q(q^2-q+1)+(p-1)q
+\sum\limits_{t=1}^{q}((q-1)\rho^{ht}+q\rho^{-2ht})\right.\\
&&\left.-\sum\limits_{t=1}^{p(q+1)}(\rho^{ht})
+\sum\limits_{t=1}^{q}(\rho^{hpt})
+\sum\limits_{t=1}^{p}(\rho^{h(q+1)t})\right]\\
&=&\frac{1}{p(q+1)}\left[q^3-q^2+q+pq-q
+q\sum\limits_{t=1}^{q}(\rho^{ht}+\rho^{-2ht})+p\right]\\
&\ge&\frac{1}{p(q+1)}\left[q^3-q^2+pq
-2q+p\right]=\frac{(q-2)q}{p}+1>0.
\end{array}
\]
Considering $\ZZ E_{\alpha,\beta}$,
for $x \in \{\alpha, \beta, \gamma\}$,
write 
$$
\mathfrak{A}_x  = \sum\limits_{t=1}^{g_x-1}((q-1)\rho^{\frac{q+1}{g_x}ht}+q\rho^{-2\frac{q+1}{g_x}ht})\equad
\mathfrak{B}_x  = \sum\limits_{t=1}^{g_x-1}(\rho^{\alpha\frac{q+1}{g_x}ht}+\rho^{\beta\frac{q+1}{g_x}ht}+\rho^{-(\alpha+\beta)\frac{q+1}{g_x}ht}).$$
Calling $M=M_{\chi_{q(q^2-q+1)}^{(h)}}(\ZZ E_{\alpha,\beta})$, we obtain
$$
M=\frac{1}{q+1}\left[q^3-q^2+q+\mathfrak{A}_\alpha+
\mathfrak{A}_\beta+\mathfrak{A}_\gamma + \mathfrak{B}_\alpha+ \mathfrak{B}_\beta+ \mathfrak{B}_\gamma -\sum\limits_{t=1}^{q}(\rho^{\alpha ht}+\rho^{\beta ht}+\rho^{-(\alpha+\beta) ht}) \right]$$
whence

$$
\begin{array}{rcl}
M&=&\frac{1}{q+1}\left[q^3-q^2+q+\sum\limits_{t=1}^{g_\alpha-1}((q-1)\theta_\alpha^{ht}+q\theta_\alpha^{-2ht})\right.\\
&&+\sum\limits_{t=1}^{g_\beta-1}((q-1)\theta_\beta^{ht}+q\theta_\gamma^{-2ht}) +\sum\limits_{t=1}^{g_\gamma-1}((q-1)\theta_\gamma^{ht}+q\theta_\gamma^{-2ht})\\
&&-\sum\limits_{t=1}^{q}(\rho^{\alpha ht}+\rho^{\beta ht}+\rho^{-(\alpha+\beta) ht})+\sum\limits_{t=1}^{g_\alpha-1}(\theta_\alpha^{\alpha ht}+\theta_\alpha^{\beta ht}+\theta_\alpha^{-(\alpha+\beta)ht})\\
&&\left.+\sum\limits_{t=1}^{g_\beta-1}(\theta_\beta^{\alpha ht}+\theta_\beta^{\beta ht}+\theta_\beta^{-(\alpha+\beta)ht})+\sum\limits_{t=1}^{g_\gamma-1}(\theta_\gamma^{\alpha ht}+\theta_\gamma^{\beta ht}+\theta_\gamma^{-(\alpha+\beta)ht})\right].
\end{array}
$$
In the same way of $M_{\chi_{q^2-q+1}^{(h)}}(\ZZ E_{\alpha,\beta})$, we have
$$
\begin{array}{rcl}
M&=&\frac{1}{q+1}\left[q^3-q^2+q+g_\alpha+g_\beta-2+\sum\limits_{t=1}^{g_\alpha-1}((q-1)\theta_\alpha^{ht}+q\theta_\alpha^{-2ht})\right.\\
&&+\sum\limits_{t=1}^{g_\beta-1}((q-1)\theta_\beta^{ht}+q\theta_\gamma^{-2ht}) +\sum\limits_{t=1}^{g_\gamma-1}((q-1)\theta_\gamma^{ht}+q\theta_\gamma^{-2ht})\\
&&-\sum\limits_{t=1}^{q}(\rho^{\alpha ht}+\rho^{\beta ht}+\rho^{-(\alpha+\beta) ht})
+2\sum\limits_{t=1}^{g_\alpha-1}(\theta_\alpha^{ ht})\\
&&\left.+2\sum\limits_{t=1}^{g_\beta-1}(\theta_\beta^{ ht})+\sum\limits_{t=1}^{g_\gamma-1}(2\theta_\gamma^{ ht}+\theta_\gamma^{-2ht})\right]\\
&=&\frac{1}{q+1}\left[q^3-q^2+q+g_\alpha+g_\beta-2+\sum\limits_{t=1}^{g_\alpha-1}((q+1)\theta_\alpha^{ht}+q\theta_\alpha^{-2ht})\right.\\
&&+\sum\limits_{t=1}^{g_\beta-1}((q+1)\theta_\beta^{ht}+q\theta_\gamma^{-2ht}) +(q+1)\sum\limits_{t=1}^{g_\gamma-1}(\theta_\gamma^{ht}+\theta_\gamma^{-2ht})\\
&&\left.-\sum\limits_{t=1}^{q}(\rho^{\alpha ht}+\rho^{\beta ht}+\rho^{-(\alpha+\beta) ht})\right].\\
\end{array}
$$
Thus,
\[
\begin{array}{rcl}
M&\ge&\frac{1}{q+1}\left[q^3-q^2+q+g_\alpha+g_\beta-2-4(q+1)-2q-3q\right]\\
&=&\frac{1}{q+1}\left[q^3-q^2-8q+g_\alpha+g_\beta-6\right]>\frac{q^3-q^2-8q-6}{q+1}
= q^2-2q-6,
\end{array}
\]
that is positive if $q>3$. 
Moreover, if $q=3$ we get $M_{\chi_{q(q^2-q+1)}^{(h)}}(\ZZ E_{1,3})\in\{5,7\}$.
In conclusion, for the element $\ZZ F$, we have 
\[
\begin{array}{rcl}
M_{\chi_{q(q^2-q+1)}^{(h)}}(\ZZ F)&=& \frac{1}{q^2-1}\left[q^3-q^2+q+\sum\limits_{t=1}^{q}((q-1)\rho^{ht}+q\rho^{-2ht})\right.\\
&&\left.+\sum\limits_{t=1}^{q^2-2} \rho^{th} -\sum\limits_{t=1}^{q}\rho^{(q-1)th}\right]\\
&=&\frac{1}{q^2-1}\left[q^3-q^2+(q-1)\sum\limits_{t=1}^{q} \rho^{-2ht} \right]\\
&\ge&\frac{1}{q^2-1}\left[q^3-q^2-q+1\right]=q-1>0.
\end{array}
\]
Hence, the character $\chi_{q(q^2-q+1)}$ is unisingular.

Continuing with the character $\chi_{(q-1)(q^2-q+1)}^{(n,m)}$, we can ignore the class $\ZZ G$ for the computation of $M_{\chi_{(q-1)(q^2-q+1)}^{(n,m)}}$.
Called $M=M_{\chi_{(q-1)(q^2-q+1)}^{(n,m)}}(\ZZ B)$,  we have

$$
M= \left\{
\begin{array}{ll}\frac{1}{4}\left[(q-1)(q^2-q+1)+(2q-1)-2\right]=\frac{q^3-2q^2+4q-4}{4}>0 & \text{ if } p=2,\\
\frac{1}{p}\left[(q-1)(q^2-q+1)-(p-1)\right]=\frac{q^3-2q^2+2q-p}{p}>0
& \text{ if } p\neq 2.
\end{array}\right.$$
For $\ZZ D_1$, calling $M=M_{\chi_{(q-1)(q^2-q+1)}^{(n,m)}}(\ZZ D_1)$, we obtain
\[
\begin{array}{rcl}
M&=&\frac{1}{p(q+1)}\left[(q-1)(q^2-q+1)+(p-1)(2q-1)\right.\\
&&+(q-1)\sum\limits_{t=1}^{q}(\rho^{t(n-2m)}+\rho^{t(m-2n)}+\rho^{t(n+m)})\\
&&-\sum\limits_{t=1}^{p(q+1)}(\rho^{t(n-2m)}+\rho^{t(m-2n)}+\rho^{t(n+m)})\\
&&+\sum\limits_{t=1}^{q}(\rho^{pt(n-2m)}+\rho^{pt(m-2n)}+\rho^{pt(n+m)})\\
&&\left.+\sum\limits_{t=1}^{p}(\rho^{(q+1)t(n-2m)}+\rho^{(q+1)t(m-2n)}+\rho^{(q+1)t(n+m)})\right]\\
&=&\frac{1}{p(q+1)}\left[q^3-2q^2+2q-1+2pq-p-2q+1\right.\\
&&+q\sum\limits_{t=1}^{q}(\rho^{t(n-2m)}+\rho^{t(m-2n)}+\rho^{t(n+m)})\\
&&\left.-\sum\limits_{t=1}^{p(q+1)}(\rho^{t(n-2m)}+\rho^{t(m-2n)}+\rho^{t(n+m)})+3p\right]\\
&\ge&\frac{1}{p(q+1)}\left[q^3-2q^2+2pq+2p-3q-3p(q+1)\right] = \frac{q^2-3q-p}{p},
\end{array}
\]
that is positive when $q>3$. 
If $q=3$, we have only $\chi_{14}^{(1,3)}$ and in this case $M=2$.
We now consider the element $\ZZ E_{\alpha,\beta}$.
For $x \in \{\alpha, \beta, \gamma\}$,
write 
$$\begin{array}{rcl}
\mathfrak{A}_x & = &\sum\limits_{t=1}^{g_x-1}(\rho^{\frac{q+1}{g_x}t(n-2m)}+\rho^{\frac{q+1}{g_x}t(m-2n)}+\rho^{\frac{q+1}{g_x}t(n+m)}),\\
\mathfrak{B}_x & =& \sum\limits_{t=1}^{g_x-1}\left(\rho^{\frac{q+1}{g_x}t(\alpha(n-m)-\beta m)}+\rho^{\frac{q+1}{g_x}t(\alpha(m-n)-\beta n)}+\rho^{\frac{q+1}{g_x}t(\beta(n-m)-\alpha m)}+\right.\\
&&\left.
\rho^{\frac{q+1}{g_x}t(\beta(m-n)-\alpha n)}+ \rho^{\frac{q+1}{g_x}t(\alpha n+\beta m)}+\rho^{\frac{q+1}{g_x}t(\alpha m+\beta n)}\right).
\end{array}$$
Calling $M=M_{\chi_{(q-1)(q^2-q+1)}^{(n,m)}}(\ZZ E_{\alpha,\beta})$, we obtain
$$
\begin{array}{rcl}
M &=&\frac{1}{q+1}\left [q^3-2q^2+2q-1 +(q-1) (\mathfrak{A}_\alpha+\mathfrak{A}_\beta+\mathfrak{A}_\gamma)
+\mathfrak{B}_\alpha +\mathfrak{B}_\beta+
\mathfrak{B}_\gamma\right.\\
&& -\sum\limits_{t=1}^{q}\left(\rho^{t(\alpha(n-m)-\beta m)}+\rho^{t(\alpha(m-n)-\beta n)}+\rho^{t(\beta(n-m)-\alpha m)}\right.\\
&&+\left.\left.
\rho^{t(\beta(m-n)-\alpha n)}+ \rho^{t(\alpha n+\beta m)}+\rho^{t(\alpha m+\beta n)}\right)\right]\\
&\ge&\frac{1}{q+1}\left[q^3-2q^2+2q-1-9(q-1)-6q-18\right]=  q^2-3q-10.
\end{array}$$
Therefore, if $q>5$ then $M>0$.
In the remaining cases, we have $M=4$ when $q=3$, and
$M=9$ when $q=4$. 
Lastly, for $\ZZ F$, we obtain

\[
\begin{array}{rcl}
M_{\chi_{(q-1)(q^2-q+1)}^{(n,m)}}(\ZZ F)&=&\frac{1}{q^2-1}\left[(q-1)(q^2-q+1)\right.\\
&&\left.+(q-1)\sum\limits_{t=1}^{q}(\rho^{t(n-2m)}+\rho^{t(m-2n)}+\rho^{t(n+m)})\right]\\
&\ge&\frac{1}{q^2-1}\left[q^3-2q^2+2q-1-3(q-1)\right] = q-2>0.
\end{array}
\]
Thus, the character $\chi_{(q-1)(q^2-q+1)}^{(n,m)}$ is unisingular.

For the character $\chi_{(q+1)(q^2-q+1)}^{(i)}$ we can reduce  to compute the value of $M_{\chi_{(q+1)(q^2-q+1)}^{(i)}}$ for $\ZZ D_1$, $\ZZ E_{\alpha,\beta}$ and $\ZZ F$.
We get
\[
\begin{array}{rcl}
M_{\chi_{(q+1)(q^2-q+1)}^{(i)}}(\ZZ D_1)&=&\frac{1}{p(q+1)}\left[(q+1)(q^2-q+1)+p-1 +(q+1)\sum\limits_{t=1}^{q}(\rho^{it})\right.\\
&&\left.+\sum\limits_{t=1}^{p(q+1)}(\rho^{it})  -\sum\limits_{t=1}^{q}(\rho^{ipt})
-\sum\limits_{t=1}^{p}(\rho^{i(q+1)t})\right]\\
&=&\frac{1}{p(q+1)}\left[q^3+1+p-1
+q\sum\limits_{t=1}^{q}(\rho^{it})-p\right] = \frac{q}{p}(q-1),\\
M_{\chi_{(q+1)(q^2-q+1)}^{(i)}}(\ZZ E_{\alpha,\beta}) &=& \frac{1}{q+1}\left[q^3+1
+(q+1)\sum\limits_{t=1}^{q}(\rho^{it})\right] = q(q-1).
\end{array}\]
For $\ZZ F$, we have
\[
\begin{array}{rcl}
M_{\chi_{(q+1)(q^2-q+1)}^{(i)}}(\ZZ F)&=&\frac{1}{q^2-1}\left[q^3+1
+(q+1)\sum\limits_{t=1}^{q}(\rho^{it})+\sum\limits_{t=1}^{q^2-2}(\sigma^{it}+\sigma^{-it})\right.\\
&&\left.-\sum\limits_{t=1}^{q}(\sigma^{(q-1)it}+\sigma^{-(q-1)it})\right].\\
\end{array}
\]
Since $\sigma$ is a primitive $(q^2-1)$-root of $1$, we obtain
$$\begin{array}{rcl}
M_{\chi_{(q+1)(q^2-q+1)}^{(i)}}(\ZZ F)& = &\frac{1}{q^2-1}\left[q^3-1+(q+1)\sum\limits_{i=1}^q \rho^{it} - \sum\limits_{t=1}^{q}(\rho^{it}+\rho^{-it})\right]\\
&\geq &\frac{q^3-1-(q-1)}{q^2-1}=q.
\end{array}
$$
Hence, the character $\chi_{(q+1)(q^2-q+1)}^{(i)}$ is unisingular.

Lastly, for $\chi_{(q-1)(q+1)^2}^{(j)}$ we can limit the study of $M_{\chi_{(q+1)(q+1)^2}^{(j)}}$ only for the classes $\ZZ B$, $\ZZ D_1$ and $\ZZ G$.
For $\ZZ B$, if $p=2$ we get
\[
\begin{array}{rcl}
M_{\chi_{(q-1)(q+1)^2}^{(j)}}(\ZZ B)&=&\frac{1}{4}\left[(q-1)(q+1)^2-(q+1) -2\right]=\frac{q^3+q^2-2q-4}{4}>0,
\end{array}
\]
and if $p\ne2$ we obtain
\[
\begin{array}{rcl}
M_{\chi_{(q-1)(q+1)^2}^{(j)}}(\ZZ B)&=&\frac{1}{p}\left[(q-1)(q+1)^2-(p-1)\right]=\frac{q^3+q^2-q-p}{p}>0.
\end{array}
\]
For $\ZZ D_1$, we have
\[
\begin{array}{rcl}
M_{\chi_{(q-1)(q+1)^2}^{(j)}}(\ZZ D_1)&=&\frac{1}{p(q+1)}\left[(q-1)(q+1)^2-(p-1)(q+1)\right]=\frac{q^2-p}{p}.
\end{array}
\]
Finally, for $\ZZ G$, we have

\[
\begin{array}{rcl}
M_{\chi_{(q-1)(q+1)^2}^{(j)}}(\ZZ G)&=&\frac{1}{q^2-q+1}\left[(q-1)(q+1)^2-\sum\limits_{t=1}^{q^2-q}(\tau^{tj}+\tau^{-qtj}+\tau^{q^2tj})\right].
\end{array}
\]
Since $\tau$ is a primitive $(q^2-q+1)$-root of $1$, also $\tau^{-q}$ and $\tau^{q^2}$ are.
Thus,
\[
\begin{array}{rcl}
M_{\chi_{(q-1)(q+1)^2}^{(j)}}(\ZZ G)&=&\frac{1}{q^2-q+1}\left[q^3+q^2-q-1+3\right]=\frac{q^3+q^2-q+2}{q^2-q+1}=q+2.
\end{array}
\]
Therefore, the character $\chi_{(q-1)(q+1)^2}^{(j)}$ is unisingular.

\subsection{The group $\PGU_3(q)$ with $\gcd(q+1,3)=3$}

A complete set of representatives for the conjugacy classes of $\PGU_3(q)$ is the same of the case $\gcd(q+1,3)=1$, with the difference that the cardinality of the conjugacy classes of $\ZZ E_{\frac{q+1}{3},\frac{2(q+1)}{3}}$ and $(\ZZ G)^{\frac{e(q+1)}{2}}$, where $e\in\{1,2\}$, is one third of the cardinality of the other conjugacy classes of type $\ZZ E_{\alpha,\beta}$ and $(\ZZ G)^d$, respectively.

Let $\mathscr{D}'$ be the quotient set $\{d:1\le d\le q^2+q \text{ and } \frac{q+1}{3}\nmid d\}/\Delta$.
From the character table of $\GU_3(q)$ described in \cite{En}, we can build the character table of $\PGU_3(q)$, that is Table~\ref{PGU3n2}, where 
$\xi=\rho^{\frac{q+1}{3}}$ is a primitive cubic root of $1$,
$k\in\{0,1,2\}$, $ e\in\{1,2\}$, $ 1\le h \le q$  with $\frac{q+1}{3}\nmid h$, $(n,m)\in\mathscr{E}$, $i\in \mathscr{C}$ and $ j\in \mathscr{D}'$.
Note that $|\mathscr{C}|=\frac{(q+1)(q-2)}{2}$, $|\mathscr{D}'|=\frac{(q+1)(q-2)}{3}$ and $|\mathscr{E}|=\frac{(q+1)(q-2)}{6}+1$.

\begin{table}[tph]
	\rotatebox{90}{
		\begin{footnotesize}
			$\begin{array}{c|cccccc}	
				& \ZZ I & \ZZ A & \ZZ B & (\ZZ C)^a & \ZZ D_b  & \ZZ E_{\frac{1}{3}(q+1),\frac{2}{3}(q+1)} \\ \hline
				
				\chi_{1}^{(k)}   & 1 & 1 & 1 & \xi^{ak} & \xi^{bk} & 1 \\\hline
				
				\chi_{q^2-q}^{(k)}   & q^2-q  & -q & 0 & -(q-1)\xi^{ak} & \xi^{bk} & 2 \\\hline
				
				\chi_{q^3}^{(k)}   & q^3 & 0 & 0 & q\xi^{ak} & 0 & -1 \\\hline
				
				\chi_{q^2-q+1}^{(h)}   & q^2-q+1 & -(q-1) & 1 & -(q-1)\rho^{ah}+\rho^{-2ah} & \rho^{bh}+\rho^{-2bh} & 3 \\\hline
				
				\chi_{q(q^2-q+1)}^{(h)}   & q(q^2-q+1) & q & 0 & (q-1)\rho^{ah}+q\rho^{-2ah} & -\rho^{bh} & -3\\\hline
				
				\chi_{(q-1)(q^2-q+1)}^{(n,m)}   & (q-1)(q^2-q+1) & 2q-1 & -1 & (q-1)(\rho^{an}+\rho^{am}+\rho^{-a(n+m)}) & 
                \begin{array}{c}
                -\rho^{bn}-\rho^{bm}-\\
                \rho^{-b(n+m)} 
                \end{array}
                 & -3\xi^{n+2m}-3\xi^{m+2n} \\\hline
				
				\chi_{(q+1)(q^2-q+1)}^{(i)}   & (q+1)(q^2-q+1) & 1 & 1 & (q+1)\rho^{-ai} & \rho^{-bi} & 0 \\\hline
				
				\chi_{(q-1)(q+1)^2}^{(j)}   & (q-1)(q+1)^2 & -(q+1) & -1 & 0 & 0 & 0 
			\end{array}$
	\end{footnotesize}}
    \quad
	\rotatebox{90}{
		\begin{footnotesize}
			$\begin{array}{c|cccc}
				& \ZZ E_{\alpha,\beta} & (\ZZ F)^c & (\ZZ G)^{\frac{e}{3}(q-1)} &  (\ZZ G)^d  \\ \hline
				
				\chi_{1}^{(k)}   & \xi^{(\alpha+\beta)k} & \xi^{ck} & \xi^{ek} & \xi^{dk}\\\hline
				
				\chi_{q^2-q}^{(k)}   & 2\xi^{(\alpha+\beta)k} & 0 & -\xi^{ek} & -\xi^{dk} \\\hline
				
				\chi_{q^3}^{(k)}   & -\xi^{(\alpha+\beta)k} & \xi^{ck} & -\xi^{ek} & -\xi^{dk} \\\hline
				
				\chi_{q^2-q+1}^{(h)}   & \rho^{(\alpha-2\beta)h}+ \rho^{(\beta-2\alpha)h} + \rho^{(\alpha+\beta)h} & \rho^{ch} & 0 & 0 \\\hline
				
				\chi_{q(q^2-q+1)}^{(h)}  & -\rho^{(\alpha-2\beta)h}- \rho^{(\beta-2\alpha)h} - \rho^{(\alpha+\beta)h} & \rho^{ch} & 0 & 0 \\\hline
				
				\chi_{(q-1)(q^2-q+1)}^{(n,m)}   & 
                \begin{array}{c}-\rho^{\beta m-\alpha(n+m)}-\rho^{\beta n-\alpha(n+m)}-\rho^{\alpha m-\beta(n+m)}-\\\rho^{\alpha n-\beta(n+m)}- \rho^{\alpha n+\beta m}-\rho^{\alpha m+\beta n}\end{array} & 0 & 0 & 0 \\\hline
				
				\chi_{(q+1)(q^2-q+1)}^{(i)}   & 0 & \sigma^{ci}+\sigma^{-qci} & 0 & 0 \\\hline
				
				\chi_{(q-1)(q+1)^2}^{(j)}   & 0 & 0 & -3\xi^{ek} & -\tau^{dj}-\tau^{-qdj}-\tau^{q^2 dj} \\
			\end{array}$
	\end{footnotesize}}
	\caption{Character table of $\PGU_3(q)$, with $\gcd(q+1,3)=3$.}
	\label{PGU3n2}
\end{table}

We can limit our study only to the characters $\chi^{(k)}_{q^2-q}$ and $\chi^{(k)}_{q^3}$ where $k\in\{1,2\}$, since $\chi_1^{(1)}$ and $\chi_1^{(2)}$ are two non principal linear characters, thus not unisingular, and all other characters are the same of the case $\gcd(q+1,3)=1$.

Since in the previous section we did not evaluate $M_\chi$ for some characters $\chi$ when $q\in \{2,5\}$ and since there are no characters of type $\chi_{q^2-q+1}^{(h)}$ for $q=2$, we must compute the value of $M_{\chi_{q^2-q+1}^{(h)}}(\ZZ E_{\alpha,\beta})$ for $q=5$, the value of  $M_{(q-1)(q^2-q+1)}^{(n,m)}(\ZZ D_1)$ for $q=2$ and the value of $M_{(q-1)(q^2-q+1)}^{(n,m)}(\ZZ E_{\alpha,\beta})$  for both $q=2$ and $q=5$. 
Then, for $q=2$ we have $M_{(q-1)(q^2-q+1)}^{(n,m)}(\ZZ D_1)=1$ and $M_{(q-1)(q^2-q+1)}^{(n,m)}(\ZZ E_{1,2})=3$. 
Instead, for $q=5$ we obtain $M_{\chi_{q^2-q+1}^{(h)}}(\ZZ E_{\alpha,\beta})\in \{3,5, 9\}$ and $M_{(q-1)(q^2-q+1)}^{(n,m)}(\ZZ E_{\alpha,\beta})\in\{13,14,16,18,30\}$.

For the character $\chi^{(k)}_{q^2-q}$ we can avoid computing $M_{\chi^{(k)}_{q^2-q}}(\ZZ B)$ when $p\ne2$.
For $p=2$ we have 
$$
M_{\chi^{(k)}_{q^2-q}}(\ZZ B)= \frac{q^2-q-q}{4} =\frac{q^2-2q}{4},$$
that is positive if and only if $q\ne 2$.
For $\ZZ D_1$, since $\xi^k$ and $\xi^{kp}$ are primitive cubic roots of $1$ and $\xi^{q+1}=1$, we get
\[
\begin{array}{rcl}
M_{\chi^{(k)}_{q^2-q}}(\ZZ D_1)&=&\frac{1}{p(q+1)}\left[q^2-q-q(p-1)-(q-1)\sum\limits_{t=1}^q(\xi^{tk})\right.\\
&&\left.+\sum\limits_{t=1}^{p(q+1)}(\xi^{tk})-\sum\limits_{t=1}^{q}(\xi^{ptk})-\sum\limits_{t=1}^{p}(\xi^{(q+1)tk})\right],\\
&=& \frac{1}{p(q+1)}\left[q^2-pq+q-1+1-p\right]=\frac{q^2-pq+q-p}{p(q+1)}=q-p,
\end{array}
\]
that is $0$ if and only if $q=p$.
For $\ZZ E_{\alpha,\beta}$, we obtain

\[
\begin{array}{rcl}
M_{\chi^{(k)}_{q^2-q}} (\ZZ E_{\alpha,\beta})&=&\frac{1}{q+1}\left[q^2-q-(q-1)\sum\limits_{t=1}^{g_\alpha-1}(\xi^{\frac{q+1}{g_\alpha}tk})
-(q-1)\sum\limits_{t=1}^{g_\beta-1}(\xi^{\frac{q+1}{g_\beta}tk})\right.\\
&& -(q-1)\sum\limits_{t=1}^{g_\gamma-1}(\xi^{\frac{q+1}{g_\gamma}tk}) +2\sum\limits_{t=1}^{q}(\xi^{t(\alpha+\beta)k})-2\sum\limits_{t=1}^{g_\alpha-1}(\xi^{\frac{q+1}{g_\alpha}t(\alpha+\beta)k})\\
&&\left.-2\sum\limits_{t=1}^{g_\beta-1}(\xi^{\frac{q+1}{g_\beta}t(\alpha+\beta)k})-2\sum\limits_{t=1}^{g_\gamma-1}(\xi^{\frac{q+1}{g_\gamma}t(\alpha+\beta)k})\right]\\
 &\ge&\frac{1}{q+1} \left[q^2-q-(q-1)(g_\alpha+g_\beta+g_\gamma-3)-2\right.\\
 &&\left.-2(g_\alpha+g_\beta+g_\gamma-3)\right]\\
 &=&\frac{1}{q+1}\left[q^2-q-2-(q+1)(g_\alpha+g_\beta+g_\gamma-3)\right].
\end{array}
\]
By Lemma~\ref{DisugualianzaSommaProdotto}, when $q>5$ we have
$$
M_{\chi^{(k)}_{q^2-q}} (\ZZ E_{\alpha,\beta})  >  \frac{1}{q+1}\left[q^2-q-2-(q+1)(q-2)\right]=0.
$$
Moreover, if $q=2$, then $M_{\chi^{(k)}_{q^2-q}} (\ZZ E_{1,2})=2$, and if $q=5$ then $M_{\chi^{(k)}_{q^2-q}} (\ZZ E_{\alpha,\beta})\in\{3,4,8\}$.
For the remaining classes, we get 
$$
\begin{array}{rcl}
M_{\chi^{(k)}_{q^2-q}} (\ZZ F) &=& \frac{1}{q^2-1}\left[q^2-q-(q-1)\sum\limits_{t=1}^q(\xi^{tk})\right] = \frac{q^2-q+q-1}{q^2-1} =1>0,\\
M_{\chi^{(k)}_{q^2-q}} (\ZZ G) &=& \frac{1}{q^2-q+1}\left[q^2-q-\sum\limits_{t=1}^{q^2-q}(\xi^{tk})\right]=\frac{q^2-q+1}{q^2-q+1}=1.
\end{array}$$
Hence $\chi^{(k)}_{q^2-q}$, with $k\in \{1,2\}$ is unisingular if and only if $q \ne p$.

For the character $\chi^{(k)}_{q^3}$, we can ignore the class $\ZZ B$ to compute the value of $M_{\chi^{(k)}_{q^3}}$.
For $\ZZ D_1$, we have
$$
M_{\chi^{(k)}_{q^3}}(\ZZ D_1) = \frac{1}{p(q+1)}\left[q^3+q\sum\limits_{t=1}^q(\xi^{tk})\right]=\frac{q^3-q}{p(q+1)}=\frac{q(q-1)}{p}.
$$
For $\ZZ E_{\alpha,\beta}$, we obtain
\[
\begin{array}{rcl}
M_{\chi^{(k)}_{q^3}}(\ZZ E_{\alpha,\beta})&=&\frac{1}{q+1}\left[q^3+q\sum\limits_{t=1}^{g_\alpha-1}(\xi^{\frac{q+1}{g_\alpha}tk})+q\sum\limits_{t=1}^{g_\beta-1}(\xi^{\frac{q+1}{g_\beta}tk})+q\sum\limits_{t=1}^{g_\gamma-1}(\xi^{\frac{q+1}{g_\gamma}tk})\right.\\
&&-\sum\limits_{t=1}^{q}(\xi^{t(\alpha+\beta)k})+\sum\limits_{t=1}^{g_\alpha-1}(\xi^{\frac{q+1}{g_\alpha}t(\alpha+\beta)k})\\
&&\left.+\sum\limits_{t=1}^{g_\beta-1}(\xi^{\frac{q+1}{g_\beta}t(\alpha+\beta)k})+\sum\limits_{t=1}^{g_\gamma-1}(\xi^{\frac{q+1}{g_\gamma}t(\alpha+\beta)k})\right]\\
&\ge&\frac{q^3-3q-q-3}{q+1} =\frac{q^3-4q-3}{q+1}=q^2-q-3,
\end{array}\]
that is positive if $q>2$.
On the other hand, if $q=2$, then $M_{\chi^{(k)}_{q^3}}(\ZZ E_{1,2})=2$. For $\ZZ F$, since $\xi^{q-1}$ is a primitive cubic root of $1$, we get

$$M_{\chi^{(k)}_{q^3}} (\ZZ F) = \frac{1}{q^2-1}\left[q^3+q\sum\limits_{t=1}^q(\xi^{tk})+\sum\limits_{t=1}^{q^2-2}(\xi^{tk})-\sum\limits_{t=1}^q(\xi^{(q-1)tk})\right] = \frac{q^3-q}{q^2-1}=q.$$
Lastly, for $\ZZ G$, we obtain
$$
M_{\chi^{(k)}_{q^3}} (\ZZ G)  = \frac{1}{q^2-q+1}\left[q^3-\sum\limits_{t=1}^{q^2-q}(\xi^{tk})\right]=\frac{q^3+1}{q^2-q+1}=q+1.$$
Thus, the character $\chi^{(k)}_{q^3}$ with $k\in \{1,2\}$ is unisingular.

We can resume the results of these two sections with the following.
\begin{prop}\label{pPGU3}
    Let $\chi$ be an irreducible character of $\PGU_3(q)$.
    Then $\chi$ is not unisingular if and only if one of the following cases occurs:
    \begin{enumerate}[$(1)$]
        \item $\chi$ is a nontrivial linear character;

        \item $\chi(1)=q^2-q$,
        where $\chi$ is rational-valued if 
        $q\equiv 2 \pmod 3$ is not a prime;

        \item $q$ is an odd prime, $\chi(1)=q^2-q+1$ and $\chi$ is  not rational-valued.

    \end{enumerate}
\end{prop}

\subsection{The group $\PSU_3(q)$}

The character table and the conjugacy classes of $\PSU_3(q)$ have been described by Simpson and Frame in \cite{SF} (see also \cite{G,O}).
Since $\PSU_3(q)=\PGU_3(q)$ when $\gcd(q+1,3)=1$, we can assume that $3 \mid (q+1)$.

The center $\ZZ$ of $\SU_3(q)$ is $\{I,\omega I,\omega^2 I\}$, where $\omega$ is an element of order $3$ in $\F_{q^2}^*$. 
In this case, let $\lambda$ be an element of order $q^2-q+1$ in $\F_{q^6}^*$ and $\eta$ be an element of order $q^2-1$ in $\F_{q^2}^*$.
Moreover, let $\nu=\eta^{q-1}$, that is an element of order $q+1$ in $\F_{q}^*$. 
Let us consider the following matrices in $\SL_3(q^6)$:
$$
	\overline{A}=\begin{pmatrix}
		1 & 0 & 0\\
		1 & 1 & 0\\
		0 & 0 & 1\\
	\end{pmatrix}, \;\;
	\overline{B_f}=\begin{pmatrix}
		1 & 0 & 0\\
		1 & 1 & 0\\
		0 & \nu^f & 1\\
	\end{pmatrix}, \;\;
	\overline{C} =\begin{pmatrix}
		\nu & 0 & 0\\
		0 & \nu & 0\\
		0 & 0 & \nu^{-2}\\
	\end{pmatrix},\;\;
	\overline{D_b} =\begin{pmatrix}
		\nu^b & 0 & 0\\
		1 & \nu^b & 0\\
		0 & 0 & \nu^{-2b}
	\end{pmatrix},$$
    $$
	\overline{E_{\alpha,\beta}}=\begin{pmatrix}
		\nu^\alpha & 0 & 0\\
		0 & \nu^\beta & 0\\
		0 & 0 & \nu^{-(\alpha+\beta)}\\
	\end{pmatrix},\;\;
	\overline{F}=\begin{pmatrix}
		\eta & 0 & 0\\
		0 & \eta^{-q} & 0 \\
		0 & 0 & \nu
	\end{pmatrix},\;\;
	\overline{G} = \begin{pmatrix}
		\lambda & 0 & 0\\
		0 & \lambda^{-q} & 0\\
		0 & 0 & \lambda^{q^2}
	\end{pmatrix},$$
    where  $0\le f\le 2$, $1\le b\le \frac{q-2}{3}$ and $1\le\alpha<\beta\le \frac{q+1}{3}$.
Let $A,B,C,D_b,E_{0,\frac{q+1}{3}},E_{\alpha,\beta}, F,G$ be elements of $\SU_3(q)$ conjugated in $\GL_3(q^6)$ respectively to $\overline{A},\overline{B}, \overline{C},\overline{D_b},\overline{E_{0,\frac{q+1}{3}}},\overline{E_{\alpha,\beta}},\overline{F}$ and $\overline{G}$.
We will call $\mathscr{C}''$ the quotient set $\left\{ c:1\le c\le \frac{q^2-1}{3},(q-1)\nmid c\right\}/\sim$ and $\mathscr{D}''$ the quotient set $\left\{d:1\le d\le \frac{q^2-q -2}{3}\right\}/\vartriangle$.
Note that $|\mathscr{C}''|=\frac{(q+1)(q-2)}{6}$, $|\mathscr{D}''|=\frac{(q+1)(q-2)}{9}$.

A complete representative set for the conjugacy classes of $\PSU_3(q)$ is given by
\[
\ZZ I, \quad \ZZ A, \quad \ZZ B_f, \quad (\ZZ C)^a, \quad \ZZ D_b,\quad \ZZ E_{\alpha,\beta},\quad \ZZ E_{0,\frac{q+1}{3}},\quad (\ZZ F)^c \equad (\ZZ G)^d,
\]
where
$0\le f\le 2$, $1\le a\le \frac{q-2}{3}$, $1\le b\le \frac{q-2}{3}$, 
$1\le\alpha<\beta\le \frac{q+1}{3}$, $c\in \mathscr{C}''$ and $d\in \mathscr{D}''$.
As in the previous sections, $\ZZ E_{\alpha,\beta}$ is always a power of some $\ZZ E_{\gamma,\delta}$, with $\gcd(\gamma,\delta)=1$.
Then we can study only $M_\chi(\ZZ E_{\alpha,\beta})$ with $\gcd(\alpha,\beta)=1$.
From now, we will suppose that $\alpha$ and $\beta$ are coprime.
Moreover, $\ZZ D_b$ is conjugate to a power $(\ZZ D_1)^k$.
Hence, we will consider only $M_\chi (\ZZ D_1)$.
The order of elements of each nontrivial conjugacy class is
\[
o(\ZZ A) = p, \quad o(\ZZ B_f) = \begin{cases}
	4 & \text{if } p=2 \\
	p & \text{otherwise}
\end{cases},\quad 
o(\ZZ C)= \frac{q+1}{3},\quad
o(\ZZ D_1)= \frac{p(q+1)}{3},
\]
\[ o(\ZZ E_{\alpha,\beta})= q+1, 
\quad 
o(\ZZ E_{0,\frac{q+1}{3}})=3,\quad
o(\ZZ F)= \frac{q^2-1}{3}, \quad o(\ZZ G) =\frac{q^2-q+1}{3}.\]

About the power maps we have: 
\begin{itemize}
	\item the element $(\ZZ A)^k$ is conjugate to $\ZZ A$ for any $1\le k\le p-1$; 
	\item if $p\ne 2$ then $(\ZZ B_f)^k$ is conjugate to $\ZZ B_f$ for any $1\le k\le p-1$; if $p=2$, then $(\ZZ B_f)^3$ is conjugate to $\ZZ B_f$, while $(\ZZ B_f)^2$ is conjugate to $\ZZ A$;
	\item if $p \mid k$, then $(\ZZ D_1)^k$ is conjugate to $(\ZZ C)^k$; if $\frac{q+1}{3} \mid k$, then $(\ZZ D_1)^k$ is conjugate to $\ZZ A$; 
	\item 
	if $k(2\alpha+\beta) \equiv 0 \pmod{q+1}$, then $(\ZZ E_{\alpha,\beta})^k$ is conjugate to $(\ZZ C)^{k\alpha}$;
	if $k(\alpha+2\beta)\equiv 0 \pmod{q+1}$, then $(\ZZ E_{\alpha,\beta})^k$ is conjugate to $(\ZZ C)^{k\beta}$; 
	if $k\alpha\equiv k\beta \pmod{q+1}$, then  $(\ZZ E_{\alpha,\beta})^k$ is conjugate to $(\ZZ C)^{k\alpha}$;
	\item the element $(\ZZ E_{0,{\frac{q+1}{3}}})^2$ is conjugate to $\ZZ E_{0,{\frac{q+1}{3}}}$;
	\item if $(q-1) \mid c$, then $(\ZZ F)^c$ is  conjugate to $(\ZZ C)^{-\frac{c}{q-1}}$.
\end{itemize}
Furthermore, there are no other intersections between powers of elements belonging to different conjugacy classes.

Since $|\PGU_3(q):\PSU_3(q)|=3$, by Clifford's Theorem, any character of $\PGU_3(q)$ restricted to $\PSU_3(q)$ is either irreducible  or splits into $3$ irreducible characters of $\PSU_3(q)$ conjugated in $\PGU_3(q)$.
Thus, we must study the restrictions to $\PSU_3(q)$ of the characters $\chi_{q^2-q}^{(k)}$ and $\chi^{(h)}_{q^2-q+1}$ of $\PGU_3(q)$ and the characters of $\PSU_3(q)$ of degree $\frac{(q-1)(q^2-q+1)}{3}$.
We observe that if $q=2$ then the non unisingular irreducible characters of $\PSU_3(q)$ are the non principal linear characters and the character of degree $2$, that is $q^2-q$.
Thus, from now, we assume $q>2$.

The restrictions of the characters $\chi_{q^2-q}^{(k)}$ of $\PGU_3(q)$ are all equal to the irreducible character of $\PSU_3(q)$ of degree $q^2-q$.
Calling $\chi_{q^2-q}$ this character, it suffices to compute the value of $M_{\chi_{q^2-q}}$ only on the element $\ZZ D_1$ when $p=q$. 
Indeed, first we denote the elements $\ZZ D_1$, $\ZZ F$ and $\ZZ G$ of $\PGU_3(q)$ described in the previous section respectively with $\dot{D}$, $\dot{F}$ and $\dot{G}$.
Note that they do not belong to $\PSU_3(q)$.
Moreover, the character $\chi_{q^2-q}^{(0)}$ fails to be unisingular on the classes $\dot{D}$, $\dot{F}$ and $\dot{G}$, but $\chi_{q^2-q}^{(1)}$ and $\chi_{q^2-q}^{(2)}$ fail only on $\dot{D}$.
Therefore, we should only compute $M_{\chi_{q^2+q}}(\ZZ D_1)$.
The required values for its computation
are
$$\chi_{q^2-q}(\ZZ A) = -q,\quad \chi_{q^2-q}(\ZZ C^a)= -(q-1),\quad \chi_{q^2-q}(\ZZ D_b) = 1.$$
Thus,
\[
\begin{array}{rcl}
	M_{\chi_{q^2-q}}(\ZZ D_1)&=&\frac{3}{p(q+1)}\left[q^2-q-(p-1)q-(q-1)\frac{q-2}{3}+\frac{p(q+1)}{3}-\frac{q-2}{3}-p \right]\\
	&=&\frac{2q^2-2pq+2q-2p}{p(q+1)}=2\frac{q-p}{p},
\end{array}
\]
that is equal to $0$ if and only if $q=p$.
Hence, the character $\chi_{q^2-q}$ is not unisingular if and only if $q\ne p$.

Also the restriction of the character $\chi^{(h)}_{q^2-q+1}$ is an irreducible character of $\PSU_3(q)$.
Since $\chi^{(h)}_{q^2-q+1}$ is not unisingular on $\dot{D}$ when $q=p$, we must compute the multiplicity of the eigenvalue $1$ on its restriction only on $\ZZ D_1$ when $q=p$.
Denoting the restriction again with $\chi^{(h)}_{q^2-q+1}$, its values on $\ZZ A$, $\ZZ C^a$ and $\ZZ D_b$ are:
\begin{itemize}
	\item $\chi_{q^2-q+1}^{(h)}(\ZZ A) = -(q-1)$,
	\item $\chi_{q^2-q+1}^{(h)}(\ZZ C^a)= -(q-1)\rho^{3ah}+\rho^{-6ah}$,
	\item $\chi_{q^2-q+1}^{(h)}(\ZZ D_b) = \rho^{3ah}+\rho^{-6ah}.$
\end{itemize}
Therefore, calling $M=M_{\chi_{(q^2-q+1)}^{(h)}}(\ZZ D_1)$, we get
\[
\begin{array}{rcl}
	M&=&\frac{3}{q(q+1)}\left[q^2-q+1-(q-1)(q-1)
	+\sum\limits_{t=1}^{\frac{q-2}{3}}(-(q-1)\rho^{3ht}+\rho^{-6ht})\right.\\
	&&\left. +\sum\limits_{t=1}^{q\frac{q+1}{3}}(\rho^{3ht}+\rho^{-6ht})
	-\sum\limits_{t=1}^{\frac{q-2}{3}}(\rho^{3hqt}+\rho^{-6hqt})-\sum\limits_{t=1}^{q}(\rho^{3h(q+1)t}+\rho^{-6h(q+1)t})\right],\\
 &=& \frac{3}{p(q+1)}\left[q-q\sum\limits_{t=1}^{\frac{q-2}{3}}(\rho^{3ht}) +\sum\limits_{t=1}^{q\frac{q+1}{3}}(\rho^{3ht}+\rho^{-6ht})-2q\right].
	\end{array}\]
Now, if $h=\frac{q+1}{6}$ then 
$$M=\frac{3}{q(q+1)}\left[q+q+\frac{q(q+1)}{3}-2q\right]=1.$$
Otherwise,
$$M=\frac{3}{q(q+1)}\left[q+q-2q\right]=0.$$
Thus, the character $\chi_{(q^2-q+1)}^{(h)}$ is unisingular if and only if $h\neq \frac{q+1}{6}$.

The remaining characters we must study are the three characters of degree $\frac{(q-1)(q^2-q+1)}{3}$.
The values of these characters are written in Table~\ref{CharactersPSU(3,q)}, where $0\le n\le 2$ and $\xi$ is a cubic root of unity.
\begin{table}[ht]
    $$    	\begin{array}{c|ccccc}
		&  \ZZ I &  \ZZ A & \ZZ B_f & (\ZZ C)^a & \ZZ D_b\\ \hline
		\chi_n & \frac{(q-1)(q^2-q+1)}{3} 
		& \frac{2q-1}{3}
        & q \delta_{f,n}-\frac{q+1}{3}  & q-1 & -1\\
	\end{array}$$
    $$\begin{array}{c|cccc}
        & \ZZ E_{0,\frac{q-1}{3}} & \ZZ E_{\alpha,\beta}  
        & (\ZZ F)^c & (\ZZ G)^d\\\hline
        \chi_n  & -\xi^{\frac{q-1}{3}}-\xi^{-\frac{q-1}{3}} 		& -\xi^{\alpha-\beta}-\xi^{\beta-\alpha} & 0 & 0
	\end{array}$$
    
	\caption{Characters of degree $\frac{(q-1)(q^2-q+1)}{3}$ in $\PSU_3(q)$.}
	\label{CharactersPSU(3,q)}
\end{table}

Let us fix $n$. We may study $M_{\chi_n}$ only on the elements $\ZZ D_b$, $\ZZ E_{\alpha,\beta}$ and $\ZZ B_f$ for $f\ne n$. 
In fact, the value of $\chi_n$ on $\ZZ E_{0,\frac{q-1}{3}}$ coincides with the value for the family $\ZZ E_{\alpha,\beta}$,
allowing $(\alpha,\beta)=(0,\frac{q-1}{3})$ and the value of $\chi_n$ on the other conjugacy classes is a nonnegative integer.
Therefore, calling $M=M_{\chi_n}(\ZZ B_f)$, we have

\[
M= \frac{1}{p}\left[\frac{(q-1)(q^2-q+1)}{3}-(p-1)\frac{q+1}{3}\right]=\frac{q^3-2q^2+3q-pq-p}{3p}>0
\]
for $p \ne 2$, while for $p=2$ we have
\[
M = \frac{1}{4}\left[\frac{(q-1)(q^2-q+1)}{3}+\frac{2q-1}{3}-2\frac{q+1}{3}\right]=\frac{(q-2)(q^2+2)}{12}
\]
that is positive if $q\ne2$. 
For $\ZZ D_1$, calling $M=M_{\chi_n}(\ZZ D_1)$, we get
\[
\begin{array}{rcl}
	M&=&\frac{3}{p(q+1)}\left[\frac{(q-1)(q^2-q+1)}{3}+(p-1)\frac{2q-1}{3}+\frac{q-2}{3}(q-1)-\frac{p(q+1)}{3}+\frac{q-2}{3}+p \right]\\
	&=&\frac{q^3-2q^2+2q-1+2pq-p-2q+1+q^2-3q+2-pq-p+q-2+3p}{p(q+1)}=\frac{q^2-2q+p}{p}>0.
\end{array}
\]
Now, set $g_\alpha=\gcd(2\alpha+\beta,q+1)$, $g_\beta=\gcd(\alpha+2\beta,q+1)$ and $g_\gamma=\gcd(\alpha-\beta,q+1)$.
Calling $M=M_{\chi_n}(\ZZ_{\alpha,\beta})$, we get
\[
\begin{array}{rcl}
	M & = & \frac{1}{q+1} \left[\frac{(q-1)(q^2-q+1)}{3} +
	\sum\limits_{t=1}^{g_\alpha-1}(q-1)+
	\sum\limits_{t=1}^{g_\beta-1}(q-1)+
	\sum\limits_{t=1}^{g_\gamma-1}(q-1)\right.\\
	& & -\sum\limits_{t=1}^{q}(\xi^{(\alpha-\beta)t}+\xi^{(\beta-\alpha)t})+ \sum\limits_{t=1}^{g_\alpha-1}(\xi^{\frac{q+1}{g_\alpha}(\alpha-\beta)t}+\xi^{\frac{q+1}{g_\alpha}(\beta-\alpha)t})\\
	& & \left. +\sum\limits_{t=1}^{g_\beta-1}(\xi^{\frac{q+1}{g_\beta}(\alpha-\beta)t}+\xi^{\frac{q+1}{g_\beta}(\beta-\alpha)t})
	+\sum\limits_{t=1}^{g_\gamma-1}(\xi^{\frac{q+1}{g_\gamma}(\alpha-\beta)t}+\xi^{\frac{q+1}{g_\gamma}(\beta-\alpha)t})\right] \\
	&\ge& \frac{1}{q+1} \left[\frac{(q+1)(q^2+q+1)}{3}-(q-1)(g_\alpha+g_\beta+g_\gamma-3)-2q-6\right].
\end{array}
\]	
Since $g_\alpha, g_\beta, g_\gamma < q+1$, we have
\[	
	M > \frac{1}{q+1} \left[\frac{(q+1)(q^2+q+1)}{3}-3q(q-1)-2q-6\right]=\frac{q^3-7q^2+5q-17}{3},
\]
that is positive if $q>5$.
In the case $q=5$, we have 
$M_{\chi_n}(\ZZ E_{0,2})=10$ and
$M_{\chi_n}(\ZZ E_{1,2})=6$.
Therefore, the character $\chi_n$ is unisingular.
Hence, we obtain the following.
\begin{prop}
	Let $\chi$ be an irreducible character of $\PSU_3(q)$, with $\gcd(3,q+1)=3$.
	Then $\chi$ is not unisingular if and only if one of the following cases occurs:
	\begin{enumerate}[$(1)$]
		\item $\chi$ is a nontrivial linear character;		
		\item $q$ is a prime and $\chi(1)=q^2-q$;
        \item $q$ is an odd prime, $\chi(1)=q^2-q+1$ and $\chi$ is not rational-valued.
        
	\end{enumerate}
\end{prop}

\section{Suzuki groups and Small Ree groups}

The Suzuki groups $\Sz(q^2)={}^2\mathrm{B}_2(q^2)$, where $q^2$ is an odd power of $2$, were constructed by Michio Suzuki in \cite{Suz} and in the very same paper one can find also their character table. 
The group $\Sz(2)$ is not simple, being isomorphic to 
$\Z_5:\Z_4$. The only nonlinear character $\chi$ of $\Sz(2)$ has degree $4$ and it is not unisingular, since $M_\chi(g)=0$ for any element $g$ of order $5$.
Hence, we can write $q^2=2^{2m+1}$ with $m\geq 1$. 

As stated in \cite[2.8]{Z90}, by  a direct analysis of the character table of $\Sz(q^2)$,
Zalesski obtained that, for every non trivial representation $\Phi$ of $\Sz(q^2)$, the matrix $\Phi(g)$ has $|g|$ distinct eigenvalues for each element $g \in \Sz(q^2)$. Hence, the following follows.

\begin{prop}
    Let $m$ be a positive integer.
    Any irreducible character of $\Sz(2^{2m+1})$ is unisingular.
\end{prop}

We now consider the small Ree groups ${}^2\mathrm{G}_2(q^2)$, where $q^2$ is an odd power of $3$.
The group ${}^2\mathrm{G}_2(3)$ is not simple, being isomorphic to 
$\PSL_2(8).3=\Aut(\PSL_2(8))$. 
Using \cite{GAP4} one can check that the only irreducible characters of ${}^2\mathrm{G}_2(3)$ which are not unisingular are
the two non rational-valued characters of degree $7$ and the three characters of degree $8$.
Hence, we can write $q^2=3^{2m+1}$ with $m\geq 1$, that is $q^2\geq 27$.

Most of the character table of ${}^2\mathrm{G}_2(q^2)$ was determined by Ward in
\cite{Ward}; for the complete table, we refer to \cite{Chevie}.
Note that $q^2 \equiv 3 \pmod{8}$ and $q=\sqrt{3}q^m$.
The odd integers $\frac{q^2-1}{2}, \frac{q^2+1}{4}, q^2-\sqrt{3}q+1$ and $q^2+\sqrt{3}q+1$ are pairwise coprime.
Let $P_3$ be a Sylow $3$-subgroup of ${}^2\mathrm{G}_2(q^2)$.
Let $X\in \ZZ(P_3)$ of order $3$, $T\in P_3\setminus\ZZ(P_3)$ of order $3$ and $Y\in P_3$ of order $9$ with $Y^3=X$.
Moreover, let $J$ be an involution of ${}^2\mathrm{G}_2(q^2)$, $R$ be an element of order $\frac{q^2-1}{2}$, $S$ be an element of order $\frac{q^2+1}{4}$, $V$ be an element of order $q^2-\sqrt{3}q+1$ and $W$ be an element of order $q^2+\sqrt{3}q+1$.
The existence of these elements was proved in \cite{Ward}.

Let $\mathscr{B}$ be the quotient set $\{b: 1 \le b\le q^2-\sqrt{3}q\}/\sim$, where $b \sim\overline{b}$ if and only if $\overline{b} \equiv 
b (\sqrt{3}q-1)^\kappa  \pmod{q^2-\sqrt{3}q+1}$, for some $\kappa\geq 0$.
Let $\mathscr{E}$ be the quotient set $\{e: 1\le e \le q^2+\sqrt{3}q\}/\vartriangle$, where $e \vartriangle \overline{e}$ if and only if 
$\overline{e} \equiv e(\sqrt{3}q+1)^\kappa \pmod{q^2+\sqrt{3}q+1}$, for some $\kappa\geq 0$.
Let $\mathscr{C}$ be the quotient set $\{ (c,d): 0\le c\le1, 1\le d\le \frac{q^2-3}{4}\}/ \diamond$ where
$(c,d) \diamond(\overline{c},\overline{d})$ if and only if either
$c=\overline{c}=0$ and $\overline{d}\equiv  d \left(\frac{\sqrt{3}q+1}{2}\right)^\kappa \pmod{\frac{q^2+1}{4}}$, for some $\kappa\geq 0$, 
or $c=\overline{c}=1$ and
$\overline{d}\equiv \pm d \pmod{\frac{q^2+1}{4}}$. 
Hence, $|\mathscr{B}|=\frac{q-\sqrt{3}q}{5}$, $|\mathscr{C}|=\frac{q^2-3}{6}$ and $|\mathscr{E}|=\frac{q+\sqrt{3}q}{5}$.

A complete set of representatives for the conjugacy classes is 
$$\left\{1,\; X,\;  T^{\pm 1},\; Y,\;YT^{\pm 1},\;
J,\;  JT^{\pm 1}, \;  (JR)^a,\; V^b,\;  J^c S^d,\;  W^e:\right.$$
$$\left.  1\le a \le \frac{q^2-3}{2},\; b\in\mathscr{B},\; (c,d)\in\mathscr{C},\;
e\in \mathscr{E}\right\}.$$
In particular, the order of $YT$ and $YT^{-1}$ is $9$, the order of $JT$ and $JT^{-1}$ is $6$, the order of $JR$ is $q^2-1$ and the order of $JS$ is $\frac{q^2+1}{2}$. 
By \cite{W}, we also get that $(YT)^{-1}$ is conjugate to $YT^{-1}$, $(JT)^4=T$, $(JT)^3=J$ and $(JS)^h=J^hS^h$.
Note that 
$$J S^d = \left\{
\begin{array}{ll}
(JS)^d & \text{ if $d$ is odd},\\
(JS)^{\frac{q^2+1}{4}+d } & \text{ if $d$ is even}
\end{array}
\right.
\equad 
S^d = \left\{
\begin{array}{ll}
(JS)^{\frac{q^2+1}{4}+d } & \text{ if $d$ is odd},\\
(JS)^d & \text{ if $d$ is even}.\\
\end{array}
\right.
$$

By \cite[Theorem 1]{Z90}, $M_\chi(g)\geq 1$ for all $\chi \in \Irr({}^2\mathrm{G}_2(q^2))$ and all $g \in {}^2\mathrm{G}_2(q^2)$ such that $\gcd(|g|,3)=1$.
Thus, for the study of the unisingularity of $\chi \in \Irr({}^2\mathrm{G}_2(q^2))$ we can compute  the value of $M_\chi$ only on the elements $Y$, $YT$ and  $JT$.

In Table \ref{CharacterTableRee2} we use the following notation:  $\iota$ is the imaginary unit,
$1\le j \le \frac{q^2-3}{2}$, $ n\in\mathscr{B}$, $(\ell,k)\in\mathscr{C}$ and $ m\in \mathscr{E}$. For the power maps, we have:
\begin{itemize}
    \item if $h\in\{3,6\}$ then $Y^h$ is conjugate to $X$ and if $3\nmid h$ then $Y^h$ is conjugate to $Y$;
    \item if $h \in \{3,6\}$ then $(YT)^h$ is conjugate to $X$;
    if $h\in\{1,4,7\}$ then $(YT)^h$ is conjugate to $YT$; if $h\in\{2,5,8\}$ then $(YT)^h$ is conjugate to $YT^{-1}$;
    \item the element $(JT)^2$ is conjugate to $T^{-1}$, the element $(JT)^3$ is conjugate to $J$, the element $(JT)^4$ is conjugate to $T$ and the element $(JT)^5$ is conjugate to $JT^{-1}$.
\end{itemize}

\begin{table}[htp]
    \begin{footnotesize}
        $\begin{array}{c|ccccc} 
		&I & X& T& T^{-1}& Y\\ \hline 
		\chi_{1}&1&1&1&1&1 \\ 
		\chi_{2}&\frac{\sqrt{3}}{6}q(q^2-1) (q^2+\sqrt{3}q+1)&-\frac{1}{6}(3q^{2}+\sqrt{3}q)&\frac{\sqrt{3}}{6}(-q+q^{2}\iota)&\frac{\sqrt{3}}{6}(-q-q^{2}\iota)&\frac{\sqrt{3}}{3}q\\		
		\chi_{3}&\frac{\sqrt{3}}{6}q(q^2-1) (q^{2}-\sqrt{3}q+1)&\frac{1}{6}(3q^{2}-\sqrt{3}q)&\frac{\sqrt{3}}{6}(-q+q^{2}\iota)&\frac{\sqrt{3}}{6}(-q-q^{2}\iota)& \frac{\sqrt{3}}{3}q\\ 
		\chi_{4}&\frac{\sqrt{3}}{6}q(q^2-1) (q^2+\sqrt{3}q+1)&-\frac{1}{6}(3q^{2}+\sqrt{3}q)&\frac{\sqrt{3}}{6}(-q-q^{2}\iota)&\frac{\sqrt{3}}{6}(-q+q^{2}\iota)&\frac{\sqrt{3}}{3}q\\ 
		\chi_{5}&\frac{\sqrt{3}}{6}q(q^2-1) (q^{2}-\sqrt{3}q+1)&\frac{1}{6}(3q^{2}-\sqrt{3}q)&\frac{\sqrt{3}}{6}(-q-q^{2}\iota)&\frac{\sqrt{3}}{6}(-q+q^{2}\iota)&\frac{\sqrt{3}}{3}q\\ 
		\chi_{6}&\frac{1}{3}\sqrt{3}q(q^4-1)&-\frac{1}{3}\sqrt{3}q&\frac{\sqrt{3}}{3}(-q+q^{2}\iota)&\frac{\sqrt{3}}{3}(-q-q^{2}\iota)&-\frac{\sqrt{3}}{3}q\\ 
		\chi_{7}&\frac{1}{3}\sqrt{3}q(q^4-1)&-\frac{1}{3}\sqrt{3}q&\frac{\sqrt{3}}{3}(-q-q^{2}\iota)&\frac{\sqrt{3}}{3}(-q+q^{2}\iota)&-\frac{\sqrt{3}}{3}q\\ 
		\chi_{8}&q^{6}&0&0&0&0\\ 
		\chi_{9}& q^4-q^2+1 &-(q^2-1)&1&1&1\\ 
		\chi_{10}&q^{2}(q^4-q^2+1)&q^{2}&0&0&0\\ 
		\chi_{11}^{( j)}&(q^2+1)(q^4-q^2+1) &1&1&1&1\\ 
		\chi_{12}^{( n)}&(q^4-1)(q^2+\sqrt{3}q+1)&-(q^2+\sqrt{3}q+1)&-(\sqrt{3}q+1)&-(\sqrt{3}q+1)&-1 \\ 
		\chi_{13}^{( \ell, k)}&(q^2-1)(q^4-q^2+1) &2q^{2}-1&-1&-1&-1\\ 

		\chi_{14}^{( m)}& (q^4-1)(q^{2}-\sqrt{3}q+1)&-(q^2-\sqrt{3}q+1)&\sqrt{3}q-1&\sqrt{3}q-1&-1 \\
            \end{array}$    
            
        $\begin{array}{c|ccccc}
	& YT &YT^{-1} & J& JT& JT^{-1}\\ \hline
	\chi_{1} &1&1&1&1&1\\ 
	\chi_{2}&-\frac{\sqrt{3}}{6}q(1+\sqrt{3}\iota) &-\frac{\sqrt{3}}{6}q(1-\sqrt{3}\iota) &-\frac{1}{2}(q^2-1)&\frac{1}{2}-\frac{1}{2}q\iota&\frac{1}{2}+\frac{1}{2}q\iota\\ 
	\chi_{3} & -\frac{\sqrt{3}}{6}q(1+\sqrt{3}\iota)
        &-\frac{\sqrt{3}}{6}q(1-\sqrt{3}\iota)&\frac{1}{2}(q^2-1)&-\frac{1}{2}+\frac{1}{2}q\iota&-\frac{1}{2}-\frac{1}{2}q\iota\\ 
	\chi_{4} &-\frac{\sqrt{3}}{6}q(1-\sqrt{3}\iota)
        &-\frac{\sqrt{3}}{6}q(1+\sqrt{3}\iota)&-\frac{1}{2}(q^2-1)&\frac{1}{2}+\frac{1}{2}q\iota&\frac{1}{2}-\frac{1}{2}q\iota\\ 
	\chi_{5} &-\frac{\sqrt{3}}{6}q(1-\sqrt{3}\iota)
        &-\frac{\sqrt{3}}{6}q(1+\sqrt{3}\iota)&\frac{1}{2}(q^2-1)&-\frac{1}{2}-\frac{1}{2}q\iota&-\frac{1}{2}+\frac{1}{2}q\iota\\ 
	\chi_{6} &\frac{\sqrt{3}}{6}q(1+\sqrt{3}\iota) &\frac{\sqrt{3}}{6}q(1-\sqrt{3}\iota)&0&0&0\\ 
	\chi_{7}&\frac{\sqrt{3}}{6}q(1-\sqrt{3}\iota) & \frac{\sqrt{3}}{6}q(1+\sqrt{3}\iota)&0&0&0\\ 
	\chi_{8}&0 &0&q^{2}&0&0\\ 
	\chi_{9} &1&1&-1&-1&-1\\ 
	\chi_{10} &0&0&-q^{2}&0&0\\ 
	\chi_{11}^{( j)} &1&1&(q^2+1)(-1)^{ j}&(-1)^{ j}&(-1)^{ j}\\ 
	\chi_{12}^{( n)} & -1&-1&0&0&0\\ 
	
	\chi_{13}^{( \ell, k )} &-1&-1&(1+2(-1)^{\ell})(q^2-1)&-(1+2(-1)^{\ell}) &-(1+2(-1)^{\ell}) \\ 
	
	\chi_{14}^{(m)} &-1&-1&0&0&0
	
    \end{array}$

    \end{footnotesize}
    \caption{Partial character table of ${}^2\mathrm{G}_2(q^2)$.}
    \label{CharacterTableRee2}
\end{table}

For any $\chi \in \Irr({}^2\mathrm{G}_2(q^2))$, we get
$$\begin{array}{rcl}
M_\chi(Y) & = & \frac{1}{9}[\chi(1)+2\chi(X)+6\chi(Y)],\\
M_\chi(YT) & =& \frac{1}{9}[\chi(1)+2\chi(X)+3\chi(YT)+3\chi(YT^{-1})],\\
M_\chi(JT) & =& \frac{1}{6}[\chi(1)+\chi(J)+\chi(T)+\chi(T^{-1})+\chi(JT)+\chi(JT^{-1})].
\end{array}$$

Skipping the principal character, we start with $\chi_2$.
We get:
$$\begin{array}{rcl}
M_{\chi_2}(Y)&=&\frac{1}{9}\left[\frac{\sqrt{3}}{6}q(q^2-1) (q^2+\sqrt{3}q+1)-\frac{1}{3}(3q^{2}+\sqrt{3}q)+2\sqrt{3}q\right]\\
&=&\frac{\sqrt{3}q}{9}(q^4+\sqrt{3} q^3-3\sqrt{3} q+9)>0,\\
M_{\chi_2}(YT)&=&\frac{1}{9}\left[\frac{\sqrt{3}q^5+3q^4-3q^2-\sqrt{3}q}{6}-\frac{1}{3}(3q^{2}+\sqrt{3}q)-\sqrt{3}q\right]\\
&=& \frac{\sqrt{3}q}{9}(q^2-3)(q^2+\sqrt{3} q+3)>0,\\
M_{\chi_2}(JT)&=& \frac{1}{6}\left[\frac{\sqrt{3}q^5+3q^4-3q^2-\sqrt{3}q}{6}-\frac{q^2-1}{2}-\frac{\sqrt{3}q}{3}+1\right]\\
& = & \frac{\sqrt{3}}{36}(q+\sqrt{3}) (q^4-2\sqrt{3}q+3)>0.
\end{array}$$
Hence, the character $\chi_2$ is unisingular.
Since $\chi_4=\overline{\chi_2}$, also $\chi_4$ is unisingular.

For the character $\chi_3$ it suffices to compute $M_{\chi_3}$ only for the elements $YT$ and $JT$. 
We obtain:
$$\begin{array}{rcl}
M_{\chi_3}(YT)&=&\frac{1}{9}\left[\frac{\sqrt{3}}{6}q(q^2-1) (q^2-\sqrt{3}q+1)+\frac{1}{3}(3q^{2}-\sqrt{3}q)-\sqrt{3}q\right]\\
& =& \frac{\sqrt{3} q}{54} (q^2-3)(q^2-\sqrt{3} q+3)>0,\\
M_{\chi_3}(JT)&=& \frac{1}{6}\left[\frac{\sqrt{3}q^5-3q^4+3q^2-\sqrt{3}q}{6}+\frac{q^2-1}{2}-\frac{\sqrt{3}q}{3}-1\right] = \frac{(\sqrt{3}q-3)(q^4+2\sqrt{3}q+3)}{36}>0.
\end{array}$$
Therefore, the character $\chi_3$ is unisingular.
Since $\chi_5=\overline{\chi_3}$, also $\chi_5$ is unisingular. 

For the character $\chi_6$ 
we get:
$$\begin{array}{rcl}
M_{\chi_6}(Y) & = & \frac{1}{9}\left[\frac{\sqrt{3}}{3}q(q^4-1)-\frac{2}{3}\sqrt{3}q-2\sqrt{3}q\right] = \frac{\sqrt{3}q}{9}(q^4-9)>0,\\
M_{\chi_6}(YT) & = & \frac{1}{9}\left[\frac{\sqrt{3}}{3}q(q^4-1)-\frac{2}{3}\sqrt{3}q+\sqrt{3}q\right]= \frac{\sqrt{3}q^5}{27} >0,\\
M_{\chi_6}(JT) & = & \frac{1}{6}\left[ \frac{\sqrt{3}}{3}q(q^4-1)-\frac{2\sqrt{3}q}{3} \right]= \frac{\sqrt{3}q}{18}(q^4-3)>0.
\end{array}
$$
Thus, the character $\chi_6$ is unisingular.
Since $\chi_7=\overline{\chi_6}$, also $\chi_7$ is unisingular.

The character $\chi_8$ is clearly unisingular, since  the values for the classes we are studying are nonnegative integers.

For the character $\chi_9$ it suffices to compute $M_{\chi_9}$ only for the elements $Y$ and $JT$:
\[
\begin{array}{rcl}
M_{\chi_9}(Y) & = & \frac{1}{9}\left[ q^4-q^2+1-2(q^2-1)+6 \right]=\frac{q^4-3q^2+9}{9}>0,\\
M_{\chi_9}(JT)&=&\frac{1}{6}\left[ q^4-q^2+1 +2 -1-2 \right]=\frac{q^4-q^2}{6}>0.
\end{array}
\]
Hence, the character $\chi_9$ is unisingular.

For the character $\chi_{10}$ we can only compute $M_{\chi_{10}}$ for $JT$. 
We get:
\[
M_{\chi_{10}}(JT) =\frac{1}{6}\left[q^6-q^4+q^2-q^2\right]=\frac{q^6-q^4}{6}>0.
\]
Thus, the character $\chi_{10}$ is unisingular.

Also for the character $\chi_{11}^{(j)}$ it suffices to compute $M_{\chi_{11}^{(j)}}$ for $JT$. 
We get
\[
\begin{array}{rcl}
M_{\chi_{11}^{(j)}}(JT)&=&\frac{1}{6}\left[(q^2+1)(q^4-q^2+1)+(q^2+1)(-1)^j+2(-1)^j+2\right]\\
&\ge&\frac{1}{6}\left[q^6+1-(q^2+1)-2+2\right]=\frac{q^6-q^2}{6}>0.
\end{array}
\]
Hence, the character $\chi_{11}^{(j)}$ is unisingular.

For the character $\chi_{12}^{(n)}$,
we get:
$$\begin{array}{rcl}
M_{\chi_{12}^{(n)}}(Y)= M_{\chi_{12}^{(n)}}(YT) &=&\frac{1}{9}\left[(q^4-1)(q^2+\sqrt{3}q+1)-2(q^2+\sqrt{3}q+1)-6\right]\\
&=& \frac{q^6+\sqrt{3}q^5+q^4-3q^2-3\sqrt{3}q-9}{9}>0,\\
M_{\chi_{12}^{(n)}}(JT)&=&\frac{1}{6}\left[(q^4-1)(q^2+\sqrt{3}q+1)-2(\sqrt{3}q+1)\right]\\
&=&\frac{q^6+\sqrt{3}q^5+q^4-q^2-3\sqrt{3}q-3}{6}>0.
\end{array}$$
Therefore, the character $\chi_{12}^{(n)}$ is unisingular.

For the character $\chi_{13}^{(\ell,k)}$,
we obtain:
\[
\begin{array}{rcl}
M_{\chi_{13}^{(\ell,k)}}(Y)=M_{\chi_{13}^{(\ell,k)}}(YT) &=&\frac{1}{9}\left[(q^2-1)(q^4-q^2+1)+2(2q^2-1)-6\right]\\
&=& \frac{q^6-2q^4+6q^2-9}{9}>0,\\
M_{\chi_{13}^{(\ell,k)}}(JT)&=&\frac{1}{6}\left[q^6-2q^4+2q^2-1+(q^2-1)(1+2(-1)^\ell)\right.\\
& &\left.-2-2(1+2(-1)^\ell)\right]\\
&=&\frac{1}{6}\left[q^6-2q^4+3q^2-6+2(q^2-3)(-1)^\ell\right]\\
&\ge&\frac{1}{6}\left[ q^6-2q^4+3q^2-6-2(q^2-3)\right]= \frac{q^2(q^2-1)^2}{6}>0.
\end{array}
\]
Hence, the character $\chi_{13}^{(\ell,k)}$ is unisingular.

Lastly, for the character $\chi_{14}^{(m)}$ it suffices to evaluate $M_{\chi_{14}^{(m)}}$ for $Y$ and $YT$. 
We get:
\[
\begin{array}{rcl}
\chi_{14}^{(m)}(Y) =\chi_{14}^{(m)}(YT) &=&\frac{1}{9}\left[(q^4-1)(q^2-\sqrt{3}q+1)-2(q^2-\sqrt{3}q+1)-6 \right]\\
&=&\frac{q^6-\sqrt{3}q^5+q^4-3q^2+3\sqrt{3}q-9}{9}>0.
\end{array}
\]
Hence, the character $\chi_{14}^{(m)}$ is unisingular.

We can resume the results of this section in the following.
\begin{prop}
    Let $m$ be a positive integer.
    Any irreducible character of ${}^2\mathrm{G}_2 (3^{2m+1})$ is unisingular.
\end{prop}

\section{Sporadic groups}

By computer calculations \cite{GAP4}, we can easily deal with the almost simple sporadic groups.
In fact, in Table~\ref{tab:sporadic} we list the almost simple sporadic groups possessing non unisingular representations and the exact occurrences of nonlinear irreducible characters which are not unisingular, highlighting the classes for which the characters fail to be unisingular.
Irreducible characters and conjugacy classes are indexed according to \cite{GAP4}.

\begin{table}[ht]
   $$\begin{array}{|c|c|c|c|}\hline
   G & i & \chi_i(1) & \text{Classes}  \\ \hline
\mathrm{M}_{11} & 2,3,4 & 10 & 11a, 11b \\\hline
\mathrm{M}_{23} & 2 &  22 & 23a, 23b\\\hline
\mathrm{HS} &  2 & 22 & 20a, 20b \\\hline
\mathrm{HS}.2 & 3,4 & 22 & 20a \\\hline
\mathrm{McL} &  2 & 22 &  15a,15b, 30a, 30b\\ \hline
\mathrm{McL}.2 & 3,4 &  22 & 15a, 30a\\ \hline
\end{array} $$
    \caption{Non usingular characters in almost simple sporadic groups.}
    \label{tab:sporadic}
\end{table}

Alternating and symmetric groups have been studied in \cite{Am1,Am}.
The group $A_6$ is the only simple alternating group such that $\Aut(A_n)\neq S_n$. Using \cite{GAP4} we can deal also with this exceptional case.
In the following we describe the unisingular representations for finite almost simple groups, whose socle is not a group of Lie type.

\begin{prop}\label{spor}
Let $G$ be a finite almost simple group such that $\soc(G)$ is either an alternating group or a sporadic group.
A nonlinear irreducible character $\chi$ of $G$  is unisingular if and only if none of the following cases occurs:
\begin{enumerate}[$(1)$]
\item $G=S_6$, $\chi(1)=5$, and the value of $\chi$ on a transposition is equal to $-1$;

\item $\PGU(2,9)\leq  G\leq\Aut(A_6)$, $\chi(1)=9$, and the value of $\chi$ on an element of order $10$ is equal to $-1$; 

\item $\soc(G)=A_8$ and $\chi(1)=14$;
\item $G=S_{10}$, $\chi(1)=42$, and the value of $\chi$ on a transposition is equal to $-14$;

\item $G=A_n$, $n$ odd, and $\chi(1)=n-1$;

\item $G=S_n$, $\chi(1)=n-1$, and the value of $\chi$ on a $n$-cycle is equalto $-1$; 

\item $G=S_n$, $n$ odd, $\chi(1)=\frac{n(n-3)}{2}$, and the value of $\chi$ on a transposition is equal to $-\frac{n^2-7n+12}{2}$;

    \item $G=\mathrm{M}_{11}$ and $\chi(1)=10$;
    \item $G=\mathrm{M}_{23}$ and $\chi(1)=22$;
    \item $\soc(G)=\mathrm{HS}$ and $\chi(1)=22$;
    \item $\soc(G)=\mathrm{McL}$ and $\chi(1)=22$.
\end{enumerate}
\end{prop}

In case (2), the characters $\chi$ fail to be unisingular on the elements of order $10$.

\section*{Acknowledgments}

The results of this paper are part of second author's work for his master's thesis.
The first  author is partially supported by INdAM-GNSAGA.
The two authors would like to  thank Alexandre Zalesski, who brought this problem to their attention.

\end{document}